%% file: Paper.tex
\documentclass[twocolumn, 10pt]{autart}    

\usepackage{graphicx}
\usepackage{amsfonts}
\usepackage{amsmath}
\usepackage{epstopdf}
\usepackage{amssymb}
\usepackage[toc, page, title]{appendix}
\usepackage[pagebackref=true,breaklinks=true,letterpaper=false,colorlinks,bookmarks=false, linkcolor=black, urlcolor=blue, citecolor=black]{hyperref}
\usepackage{makeidx}  
\usepackage{multicol} 
\usepackage{graphicx}
\usepackage{subfigure}
\usepackage{array}
\newcommand{\vect}[1]{#1}
\newcommand{\mtrx}[1]{#1}
\usepackage{cancel}
\usepackage{pgfplots}

\usepackage{algorithm}
\usepackage{algorithmicx}
\usepackage[noend]{algpseudocode}

\newtheorem{lemma}{Lemma}[section]

\newtheorem{theorem}[lemma]{Theorem}

\newcommand{\prox}{\mbox{prox}}
\newcommand{\proj}{\mbox{proj}}
\newcommand{\chol}{\mbox{chol}}

\begin{document}

\begin{frontmatter}

\title{Fast Robust Methods for Singular State-Space Models} 

\author[disclaimer,UWmath]{Jonathan Jonker}\ead{jonkerjo@uw.edu},    
\author[UWamath]{Aleksandr Aravkin}\ead{saravkin@uw.edu},               
\author[UWmath]{James Burke}\ead{jvburke@uw.edu},  
\author[Padua]{Gianluigi Pillonetto}\ead{giapi@dei.unipd.it},
\author[UWAPL]{Sarah Webster}\ead{swebster@apl.washington.edu}.

\address[UWmath]{Department of Mathematics, University of Washington}  
\address[UWamath]{Department of Applied Mathematics, University of Washington}             
\address[UWAPL]{Applied Physics Lab, University of Washington}        
\address[Padua]{Department of Information Engineering, University of Padova}        

\thanks[disclaimer]{The views, opinions and/or findings expressed are those of the authors and should not be interpreted as representing the official views or policies of the Department of Defense or the U.S. Government.}


\begin{abstract}
State-space models are used in a wide range of time series analysis applications. Kalman filtering and smoothing are work-horse algorithms 
in these settings. While classic algorithms assume Gaussian errors to simplify estimation,  
 recent advances use a broad range of optimization formulations to allow outlier-robust estimation, 
as well as constraints to capture prior information.  \\
Here we develop methods on state-space models where either transition or error covariances may be singular. These models frequently arise in navigation (e.g. for `colored noise' models or deterministic integrals) and are ubiquitous in auto-correlated time series models such as ARMA. We reformulate all state-space models (singular as well as nonsingluar) as constrained convex optimization problems,  and develop an efficient algorithm for this reformulation. 
The convergence rate is {\it locally linear}, with constants that do not depend on the conditioning of the problem.\\
Numerical comparisons show that the new approach outperforms competing approaches for {\it nonsingular} models, including state of the art 
interior point (IP) methods.  IP methods converge at superlinear rates; we expect them to dominate. 
However, the steep rate of the proposed approach (independent of problem conditioning) combined with cheap iterations 
wins against IP in a run-time comparison.  This suggests that the proposed approach  
can be a {\it default choice} for estimating state space models outside of the Gaussian context for singular and nonsingular models.  
To highlight the capabilities of the new framework, we focus on 
navigation applications that use singular process covariance models, 
and analyze data from a drifting mooring as a proxy for an autonomous underwater vehicle.

\end{abstract}

\end{frontmatter}

\input{intro}

\input{formulation}

\input{vehicleModel_paper}

\input{experiment}

\bibliographystyle{abbrv}
\bibliography{kalmanSurvey,uwNav}

\input{appendix}

\end{document}

%% file: intro.tex
\section{Introduction}
\label{sec:intro}

\vspace{-.1in}

The linear state space model is widely used in tracking and navigation~\cite{ybarshalom-2001a}, control~\cite{anderson2007optimal},  signal processing \cite{AndersonMoore}, and  other time series~\cite{hyndman2002state,tsay2005analysis}. 
The model assumes linear relationships between latent states with noisy observations: 
\begin{equation} \label{eq:statespace}
\begin{aligned}
x_1 & =  x_0 + w_1
\\
x_k & =  G_kx_{k-1} + w_k, \quad  k = 2, \dots, N
\\
y_k & =  H_kx_k + v_k, \quad  k = 1, \dots, N,
\end{aligned}
\end{equation}
where $x_0$ is a given initial state estimate, $x_1, \dots, x_N$ are unknown latent states with known linear process models $G_k$, and $y_1, \dots, y_N$ are observations obtained using known linear models $H_k$. 
Data must be in the range of $H_k$; so we assume $H_k$ are surjective. \\
The errors $w_k$ and $v_k$  are assumed to be mutually independent random variables with known covariances $Q_k$ and $R_k$.
In tracking and navigation, the end goal is the estimation of the latent states $\{x_k\}$. In autocorrelated time series models (e.g. Holt-Winters c.f.~\cite{hyndman2002state}, ARMA c.f.~\cite{tsay2005analysis}), estimating the state is a necessary step to  
estimating additional parameters on which $G_k$, $H_k$, $Q_k$ and $R_k$ may depend. In both settings, estimating the state sequence $\{x_k\}$ efficiently is essential.  \\
{\bf Singular Covariances}. We are particularly interested in models where $Q_k$ and $R_k$ may be singular. These models arise in all settings where state-space formulations are used. 
In navigation, the simplest example is the DC motor~\cite[pp. 95-97]{Ljung:99}:
\begin{equation}
  \label{eq:DCmotor}
\begin{aligned}
x_{k+1} &= \left(\begin{array}{cc} 0.7 & 0 \\0.084 & 1\end{array}\right) x_k
+\left(\begin{array}{c}11.81 \\0.62 \end{array}\right) (c_k + d_k)\\
y_k &= \left(\begin{array}{cc}0 & 1\end{array}\right) x_k + v_k.
\end{aligned}
\end{equation}
Here,  $y_k$ are noisy samples of the angle of the motor shaft, $c_k$ are known inputs, and  $d_k$ denote random process disturbances. The covariance matrix $Q_k$ associated to $w_k$ has dimension 2 and rank 1.  This example is general in the sense that singular models appear any time a single source of error is integrated 
into multiple states; a pervasive phenomenon in navigation models~\cite{AndersonMoore}. \\
\begin{figure}
   \begin{tabular}{l}
   \hspace{-.3in}
 \includegraphics[scale=0.6]{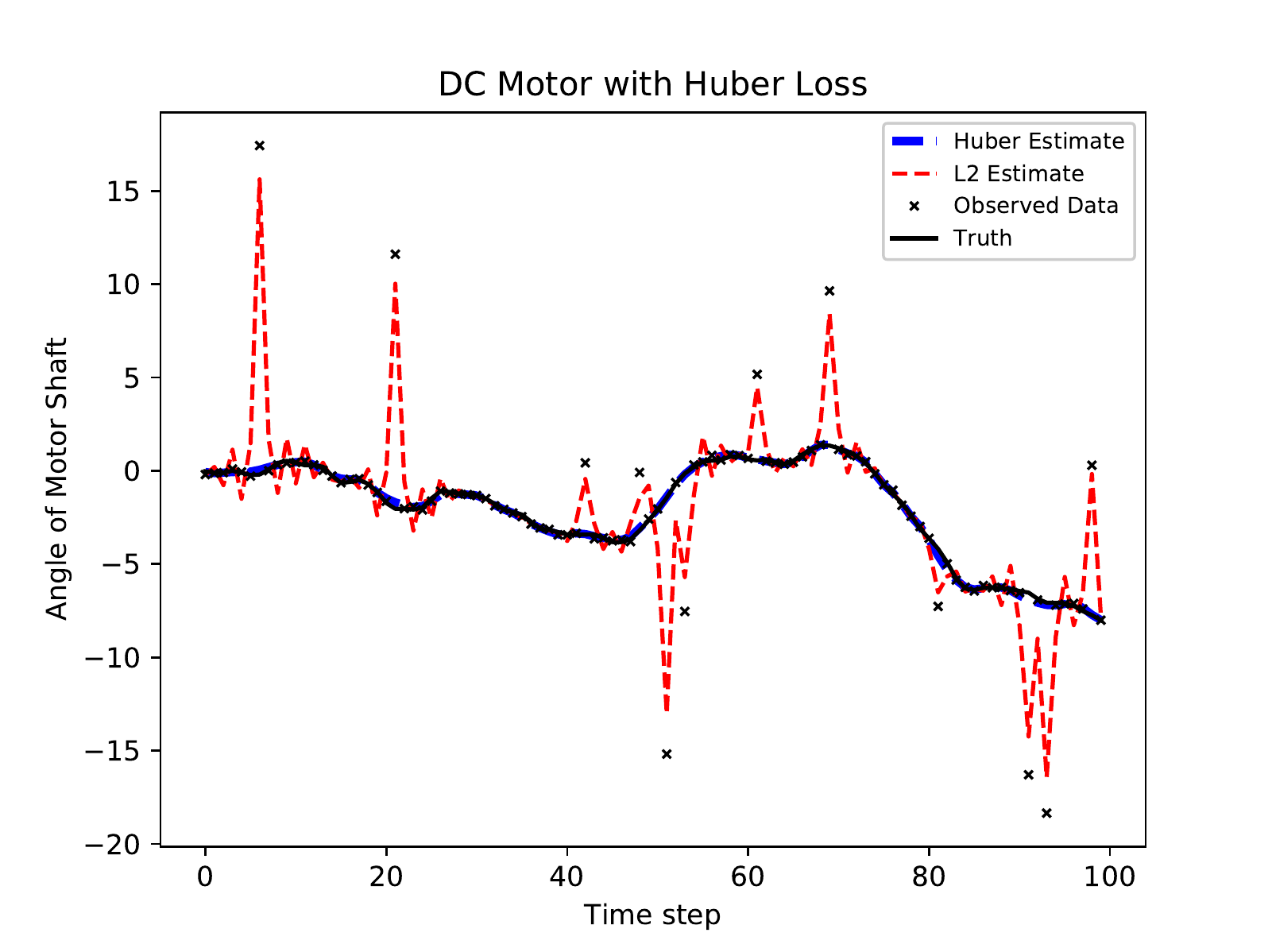} 
    \end{tabular}
 \caption{\label{fig:DCmotor} DC motor~\eqref{eq:DCmotor} with outliers, generated from a Gaussian with high variance. The process covariance $Q$ is singular, but the standard RTS smoother still finds the linear minimum variance estimate (red). Our reformulation allows using robust penalties (in this case, Huber) with a singular covariances to obtain a better solution (blue).} 
\end{figure}
The classic Kalman filter~\cite{kalman} and RTS smoother~\cite{RTS} assume that $w_k, v_k$ are Gaussian, and 
 find the {minimum variance} estimates of the state, conditioned on the observations~\cite{AndersonMoore}. 
More generally, the RTS smoother finds the linear minimum variance estimator.  This procedure is well defined for singular covariances $Q_k$ and $R_k$, and the smoother can be derived as a sequence of least squares projections~\cite{Ansley}. 
However, when the noise is not Gaussian (e.g. in the presence of outliers), these estimates are not satisfactory;
and far better estimates can be obtained through a maximum a posteriori (MAP) estimator~\cite{aravkin2017generalized}.
The results in Figure~\ref{fig:DCmotor} are obtained using the Huber loss, which is a convex penalty function that is quadratic near the origin, but with linear tails: 
\[
\rho(x) = \begin{cases}
\frac{1}{2}x^2  & \quad \text{ if } |x| \leq \kappa \\
\kappa(|x| - \frac{1}{2}\kappa) & \quad \text{ if } |x| \geq \kappa 
\end{cases}
\]
Implementing a general MAP estimator for singular covariances requires a new approach.

{\bf General Kalman Smoothing}. Classic Gaussian formulations fail when outliers are present in the data, are unable to track abrupt state changes, and cannot incorporate side information through constraints. To develop effective approaches in these cases, generalized Kalman smoothing formulations have been proposed in the last few years, see~\cite{aravkin2017generalized} and the references within. The conditional mean is no longer tractable to compute these estimates, and {\it maximum likelihood} (ML) formulations are much more natural.  The general form of Kalman smoothing considered in~\cite{aravkin2017generalized} is given by 
\begin{equation}
\label{eq:genSmoother}
\min_{x\in X} \sum_{i=0}^n \rho_1(Q_k^{-1/2}(x_k - G_kx_{k-1}))+ \rho_2(R_k^{-1/2}(y_k - H_k x_k)),
\end{equation}
where $\rho_1, \rho_2$ are convex penalties, and $x \in X$ is a set of state-space constraints. 
The two approaches agree in the nonsingular Gaussian case, where~\eqref{eq:genSmoother} becomes 
a least squares (LS) problem that can be solved with classic RTS or Mayne-Fraser smoothing algorithms~\cite{aravkin2017generalized}.

{\bf Contribution}.  We develop a new reformulation to extend~\eqref{eq:genSmoother} to {singular covariance models} $Q_k$ and $R_k$, and implement a Douglas-Rachford splitting (DRS) algorithm to solve this reformulation. The result in Figure~\ref{fig:DCmotor} uses Huber penalties for process and measurement, with the singular process covariance model from~\eqref{eq:DCmotor}.\\
We analyze the DRS for the singular reformulation, and show that it converges locally linearly for any piecewise linear quadratic (PLQ) loss, and that the rate does not depend on the conditioning of the system. {Even when the model is nonsingular, the new approach is potentially much faster} than first-order and second-order methods for~\eqref{eq:genSmoother}. The advantage increases as the models become more ill-conditioned; however the {\it local} 
linear rate means that initialization becomes very important. \\
The paper proceeds as follows. In Section~\ref{sec:survey} we discuss prior approaches to singular models. 
In Section~\ref{sec:reformulation}, we develop a constrained reformulation 
of~\eqref{eq:genSmoother}, building on early work of~\cite{Paige79} for singular least squares. 
In Section~\ref{sec:Opt}, we show how to efficiently optimize a wide range of 
singular smoothing problems using DRS. The algorithm we use has a {\it local linear rate of convergence} for any piecewise 
linear-quadratic penalties $\rho_1, \rho_2$ in~\eqref{eq:genSmoother}, and each iteration is efficiently and stably computed by exploiting dynamic problem structure. We compare the new algorithm to first-order methods, L-BFGS, and IPsolve, a toolbox specifically developed for PLQ Kalman smoothing (for nonsingular formulations). 
In Section~\ref{sec:vehicleModel}, we present a navigation model that uses singular errors. In Section~\ref{sec:numerics} we apply 
the methodology to analyze data from a drifting mooring as a proxy for an
autonomous underwater vehicle. 

\section{Related Work}
\label{sec:survey}

\vspace{-.1in}

Several approaches in the literature deal with singular models. We give a brief description and 
references for each. To ground the discussion, 
consider tracking a particle moving along a smooth path in space, where  state comprises velocity and position. 
Singular models arise naturally in this situation. 
We can model velocity as subject to error, and position as a deterministic integral: 
\begin{equation}
\label{eq:integral}
\begin{aligned}
x_{k+1} &= x_{k} + \Delta t \dot x_{k}  \\
\dot x_{k+1} &= \dot x_{k} + \epsilon_{k}.
\end{aligned}
\end{equation}
Here, the process covariance matrix $Q_k$ has rank one.\\
{\bf Using the original Kalman filter.} In the linear Gaussian setting, the original Kalman 
filter does not require $Q$ and $R$ to be invertible. Applying the Kalman filter (and RTS smoother)
will return the minimum variance estimate for singular innovation/measurement errors~\cite{AndersonMoore}. 
The limitation is that we cannot consider the general optimization context~\eqref{eq:genSmoother}, which we need 
to incorporate robustness to outliers and constraints for prior information (see example in Figure~\ref{fig:DCmotor}). \\
{\bf Changing the model.} A common approach is to modify the model to make $Q_k, R_k$ nonsingular. Treating~\eqref{eq:integral}
as a discretization of a stochastic differential equation (SDE), many authors opt  for a nonsingular error model~\cite{Jaz,Oks,Bell2008,YAA}
\[
Q_k = \begin{bmatrix} \Delta t_k & \Delta t_k^2/2 \\ \Delta t_k^2/2 & \Delta t_k^3/3\end{bmatrix},
\]
derived by computing the variance of a discretized process noise term, similar to what is done in Section~\ref{sec:vehicleModel}, see~\eqref{eq:Qderiv}. 
The approach has limitations for navigation models with high-dimensional states driven by low-dimensional errors. 
The low-dimensional error structure should simplify estimation, but instead this approach introduces full-dimensional 
and ill-conditioned $Q_k$. In addition, making $Q_k$ nonsingular is antithetical to state-space formulations for models such as ARMA, which use singularity to enforce auto-regressive constraints.   \\
{\bf Change of coordinates.}
When only $R_k$ are singular, \cite{AndersonMoore} suggests making a change of coordinates in the measurement variables and then projecting to remove the extra dimensions. The projections can vary between time points, and the approach does not extend to the singular state equation~\eqref{eq:integral}.\\
{\bf Pseudo-inverse with orthogonality constraints.}
The formulation that is closest to ours is that of~\cite{OGLB}, who replace the inverse of $Q_k$ by  a pseudo-inverse, and add orthogonality constraints (namely that projection onto the null space of $Q_k$ is zero). With potentially singular $Q_k$ and $R_k$, the maximum likelihood estimate for the Gaussian/LS model can be formulated as 
\begin{equation}
\label{eq:singularOne}
\begin{aligned}
\min_x &\sum_k ||Q_k^{\dagger/2}(x_k - G_kx_{k-1})||^2 + ||R_k^{\dagger/2}(y_k - H_kx_k)||^2\\
\quad  &\text{s.t. } Q_k^\perp(x_k - G_kx_{k-1}) = 0, \quad R_k^\perp(y_k - H_kx_k) = 0 \\
\quad & \quad \text{ for all } k = 1, \dots, N,
\end{aligned}
\end{equation}
see~\cite[Appendix A]{aravkin2017generalized}. This requires computing both the pseudo-inverse and orthogonality constraints. \\
{\bf Constrained reformulation.}
The reformulation we choose was first used by Paige \cite{Paige79}. Given the singular least squares problem 
\[
\min_{x} \|Q^{\dagger/2}(Ax-b)\|^2 \quad \mbox{s.t.} \quad Q^\perp (Ax-b) = 0,
\]
we can instead write it as 
\begin{equation}
\label{eq:singularTwo}
\min_{x, u} \|u\|^2 \quad \mbox{s.t.} \quad  \mbox Q^{1/2} u = Ax-b.
\end{equation}
It is easy to see~\eqref{eq:singularOne} and~\eqref{eq:singularTwo} are equivalent; the latter is more elegant, and only requires computing a root of $Q$, 
rather than using both $Q$ and $Q^\dagger$. 
When $Q$ is invertible, we can eliminate $u$ from both formulations 
and reduce to a least squares problem in $x$. 
Splitting the affine constraint from the original penalty has theoretical and practical advantages for general Kalman smoothing, 
as shown in the next sections.  




%% file: formulation.tex
\section{General Singular Kalman Smoothing}
\label{sec:reformulation}

\vspace{-.1in}

Following the ideas proposed by~\cite{OGLB}, we introduce variables $u_k$ for the normalized process innovations, and $t_k$ for the normalized residuals. We also introduce a penalty $\rho_3$ for the states. In the examples we consider, 
$\rho_3$ is an indicator function for the known feasible regions $X_k$:
\[
\rho_3(x_k) = \begin{cases} 0 & x_k \in X_k \\ \infty & x_k \not \in X_k\end{cases}.
\]
The reformulated singular Kalman smoothing problem is given by 
\begin{equation}
\label{eq:genKalman}
\begin{aligned}
\min_{u,t,x}  & \sum_{k=1}^N \rho_1(u_k) + \rho_2(t_k) + \rho_3(x_k)\\
 &\text{s.t.}  \quad 
 \begin{aligned}
 Q_k^{1/2}u_k &= G_kx_{k-1} - x_k\\
R_k^{1/2}t_k &= y_k - H_kx_k
\end{aligned}
\end{aligned}.
\end{equation}
This problem is equivalent to~\eqref{eq:genSmoother} when $Q_k$ and $R_k$ are nonsingular. 
For singular models,~\eqref{eq:genKalman} requires only that roots $Q^{1/2}$ and $R^{1/2}$ are available.  

{\bf Constrained Robust DC motor}. Recall the DC motor example in the introduction~\eqref{eq:DCmotor}. 
The data used to make Figure~\ref{fig:DCmotor} is contaminated with outliers, so we want to use the robust Huber loss
for the measurement errors. Suppose we also know upper and lower bounds on the states, $B := \{x: l \leq x \leq u\}$.
Then the formulation of the robust constrained singular DC motor is given by 
\[
\begin{aligned}
\min_{u,t,x} & \sum_{k=1}^N \|u_k\|^2 + \rho_h(t_k) + \delta_{B}(x_k), \quad \sigma t_k = a_k - x_{2,k},\\
  & 
 \begin{aligned}
 \begin{bmatrix}11.8 & 0 \\ 0.62 & 0 \end{bmatrix}  u_k &= x_{k+1} - \left(\begin{array}{cc} 0.7 & 0 \\0.084 & 1\end{array}\right) x_k
-\left(\begin{array}{c}11.81 \\0.62 \end{array}\right) c_k\\
\end{aligned}
\end{aligned}.
\]

{\bf Structure-preserving Reformulation.}
We now rewrite~\eqref{eq:genKalman} into a more compact form.
Define 
\begin{equation}
\label{eq:DjBj}
\begin{aligned}
D_i &= \begin{pmatrix}Q_i^{1/2} & 0 & I\\ 0 & R_i^{1/2} & H_i\end{pmatrix} \text{ for } i = 1, \dots N,\\
B_j &= \begin{pmatrix}0 & \qquad 0 & -G_{j+1}\\ 0 & \qquad 0 &  0\end{pmatrix}, \text{ for } j = 1, \dots, N-1,
\end{aligned}
\end{equation}
 and let
\begin{equation}
\label{eq:A}
A = 
\begin{pmatrix} D_1 & 0 & \dots & 0 \\ 
B_1 & D_2 & 0 & \vdots \\ 
 0& \ddots & \ddots &  0\\
  0 & 0 & B_{N-1} & D_N\end{pmatrix}.
\end{equation}
Define also 
\begin{equation}
\label{eq:zw}
\begin{aligned}
z^T &= \begin{pmatrix} u_1^T & t_1^T & x_1^T & \dots u_N^T & t_N^T & x_N^T\end{pmatrix} \\
\hat{w}^T &= \begin{pmatrix} x_0^T & y_1^T & 0 & y_2^T & \dots & 0 & y_N^T\end{pmatrix}.
\end{aligned}
\end{equation}
Now we can write~\eqref{eq:genKalman} compactly as 
\begin{equation}
\label{eq:full}
\begin{aligned}
\min_{z}  &\quad \rho(z)  \quad \text{s.t. } Az = \hat{w}, \\
 \rho(z)  &= \sum_{k=1}^N \rho_1(u_k) + \rho_2(t_k) + \rho_3(x_k).
\end{aligned}
\end{equation}
The order of blocks in $z$ is chosen to the constraint matrix $A$ in~\eqref{eq:A} lower block bi-diagonal. \\
The constraint $Az = \hat w$ raises a natural question: when is a singular Kalman smoothing model solvable? 
Clearly we want $\hat w \in\mathrm{Ran}(A)$, but we want this condition to hold for any realization of the data $\hat w$,
so we want to know when $A$ is surjective. We can characterize this condition precisely 
in terms of a simple conditions on the individual blocks $R_i, Q_i, H_i$. 

\begin{theorem}[Surjectivity of $A$] 
\label{thm:surjectivity}
The following are equivalent. 
\begin{enumerate}
\item $A$ is surjective. 
\item Each block $D_i$ is surjective. 
\item $\mathrm{null}
\left(\begin{bmatrix} Q_i^{1/2} & 0 \\ 0& R_i^{1/2} \end{bmatrix}\right) 
\subset \mathrm{Ran} \left(\begin{bmatrix} I \\ H_i\end{bmatrix}\right)$ for all $i$.
\item $R_i + H_i\left(I-(Q_i+I)^{-1}\right)H_i^T$ is invertible for all $i$. 
\end{enumerate}
\end{theorem}
The proof is given in the Appendix.

\section{Douglas-Rachford Splitting for General Singular Kalman Smoothing}
\label{sec:Opt}

\vspace{-.1in}

Consider problem~\eqref{eq:full} as a sum of two functions, $\rho +g$, with $\rho$ as in~\eqref{eq:full}
and $g$ the indicator function of the affine constraint $Az = \hat w$: 
\begin{equation}
\label{eq:indicator}
g(z) = \begin{cases} 0 & Az = \hat w\\ \infty & Az \not = \hat w\end{cases}.
\end{equation}
Douglas-Rachford splitting (DRS) is a classic algorithm for this problem. 
For a convex function $f$, define the proximity operator (see e.g.~\cite{combettes2011proximal}) as
\[
\prox_{\alpha f}(\zeta) = \arg\min_x \frac{1}{2\alpha} \|\zeta - x\|^2 + f(x).
\]
The DRS algorithm for~\eqref{eq:full} detailed in Algorithm~\ref{DRS}.
For more on splitting methods and their convergence rates see the survey~\cite{davis2016convergence}.
    \begin{algorithm}[H]
  \caption[Caption]{Douglas-Rachford Splitting (DRS)
    \label{DRS}}
  \begin{algorithmic}[1]
    \Require{Initialize at any $z^0$, $\zeta^0$.}
    \Loop
    \State {$z^k = \prox_{\tau g}(z^{k-1}-\tau \zeta^{k-1})$} 
     \State {$\zeta^k = \prox_{\sigma \rho^*}(\zeta^{k-1}+\sigma(2z^{k} - z^{k-1}))$}  
       \EndLoop
       \Return{$z^k$}
  \end{algorithmic}
\end{algorithm}
Implementing DRS in our case requires computing two proximity operators at each iteration. 
One proximity operator is $\prox_{\rho^*}$, where $\rho^*$ denotes the \textit{convex conjugate}: 
\[
\rho^*(y) = \sup_x \langle y,x \rangle - \rho(x)
\]
The prox of of a function is related to the prox of its conjugate by Moreau's decomposition:
\[
\prox_\rho(x) + \prox_{\rho^*}(x) =x.
\]
Thus it suffices to compute $\prox_\rho$.
The function $\rho$ captures all user-supplied models, including losses used process and measurement 
transitions, as well as penalties or constraints on the state, $\rho_1, \rho_2$ and $\rho_3$.
The proximity operators of these individual elements must be provided; then $\prox_\rho$ 
is a stack of these input functions.  
Proximity operators for many common functions are easily available~\cite{combettes2011proximal}, and we include a small library with our implementation\footnote{https://github.com/UW-AMO/KalmanJulia.}. \\
The second proximity operator is $\prox_g$, which is \underline{independent} of user choice for process, measurement, and prior models:
\[
\prox_g(\eta) = \arg\min_{Az = \hat w} \frac{1}{2}\|\eta-z\|^2. 
\]
This is a simple quadratic with affine constraints, with optimality conditions given by 
\[
\begin{bmatrix}I & A^T\\A & 0\end{bmatrix}\begin{bmatrix}z\\ \nu\end{bmatrix} = \begin{bmatrix}\eta\\ \hat{w}\end{bmatrix}.
\]
There are many ways to solve this system. We opt to reduce the problem to solving a block tridiagonal system:
\[
\begin{bmatrix}I & A^T\\0 & AA^T\end{bmatrix} \begin{bmatrix}z\\ \nu\end{bmatrix} = \begin{bmatrix}\eta\\ A\eta - \hat{w}\end{bmatrix}
\]
We solve  $AA^T \nu = A\eta - \hat{w}$, then back-substitute to get the optimal $z$. The system $AA^T$ does not change 
over iterations; only the right hand side changes. We can therefore compute a single factorization, then use it in each iteration. 
Since $A$ is block bidiagonal~\eqref{eq:A}, $AA^T$ is block tridiagonal; 
when $A$ is surjective, $AA^T$ is nonsingular, and 
we can find a lower block diagonal Cholesky factorization $L$  with $LL^T = AA^T$:
\begin{equation}
\label{eq:Structure}
AA^T = 
\begin{bmatrix}
a_1& b_1^T & & \\
b_1 & a_2 & b_2^T & \\
& b_2 & a_3 & b_3^T \\
&&b_3 & a_4
\end{bmatrix}, \quad L = 
\begin{bmatrix}
c_1& & & \\
d_1 & c_2 &  & \\
& d_2 & c_3 &  \\
&&d_3 & c_4
\end{bmatrix}
\end{equation}
 The factorization is detailed in Algorithm~\ref{chol}. 
\begin{algorithm}[H]
  \caption[Caption]{Block bi-diagonal Cholesky factorization for a block tri-diagonal positive definite matrix
    \label{chol}}
  \begin{algorithmic}[1]
    \Require{Input block diagonals $\{a_i\}$ and lower off-diagonals $\{b_i\}$ of block tridiagonal matrix $AA^T$~\eqref{eq:Structure}.}
    \State{$s_0 = 0, b_0 = 0$}
        \Loop{ $k=1, \dots, N$}
        \State{$s_k = a_k - b_{k-1}s_{k-1}^{-1}b_{k-1}^T$}
    \State {$c_k = \chol(s_k)$} 
     \State {$d_k = b_1 c_k^{-^T}$} 
       \EndLoop
       \Return{Diagonal blocks $\{c_i\}$ and lower-diagonal blocks $d_i$ of block $L$~\eqref{eq:Structure}}
  \end{algorithmic}
\end{algorithm}
Algorithm~\ref{chol} is derived as follows. Multiplying out $LL^T$ we have
\[
a_1 = c_1c_1^T, \quad d_1 = b_1 c_1^{-T}
\]
To compute $c_1$ we need the standard the Cholesky factorization of $a_1$. Then
\[
c_2c_2^T = a_2 - b_1a_1^{-1}b_1^T, \quad d_2 = b_2c_2^{-1}. 
\]
For convenience, we introduce the recursively defined auxiliary terms $s_k$, with $s_1 = a_1$, and
\[
s_k = a_k - b_{k-1}s_{k-1}^{-1}b_{k-1}^T.
\]
Then each $c_k$ is the standard Cholesky factorization of $s_k$, and $d_k$ is immediately computed as in Algorithm~\ref{chol}.
The overall complexity required for the single factorization is $O(n^3N)$.
Once $L$ has been pre-computed, we need only $O(n^2N)$ arithmetic operations to solve $LL^T \nu = Ac - \hat{w}$ 
for any right hand side. This is the same complexity as that of a matrix-vector multiply with $A$.


%
%
%

\begin{figure}
  \begin{center}
   \begin{tabular}{c}
 \includegraphics[scale=0.5]{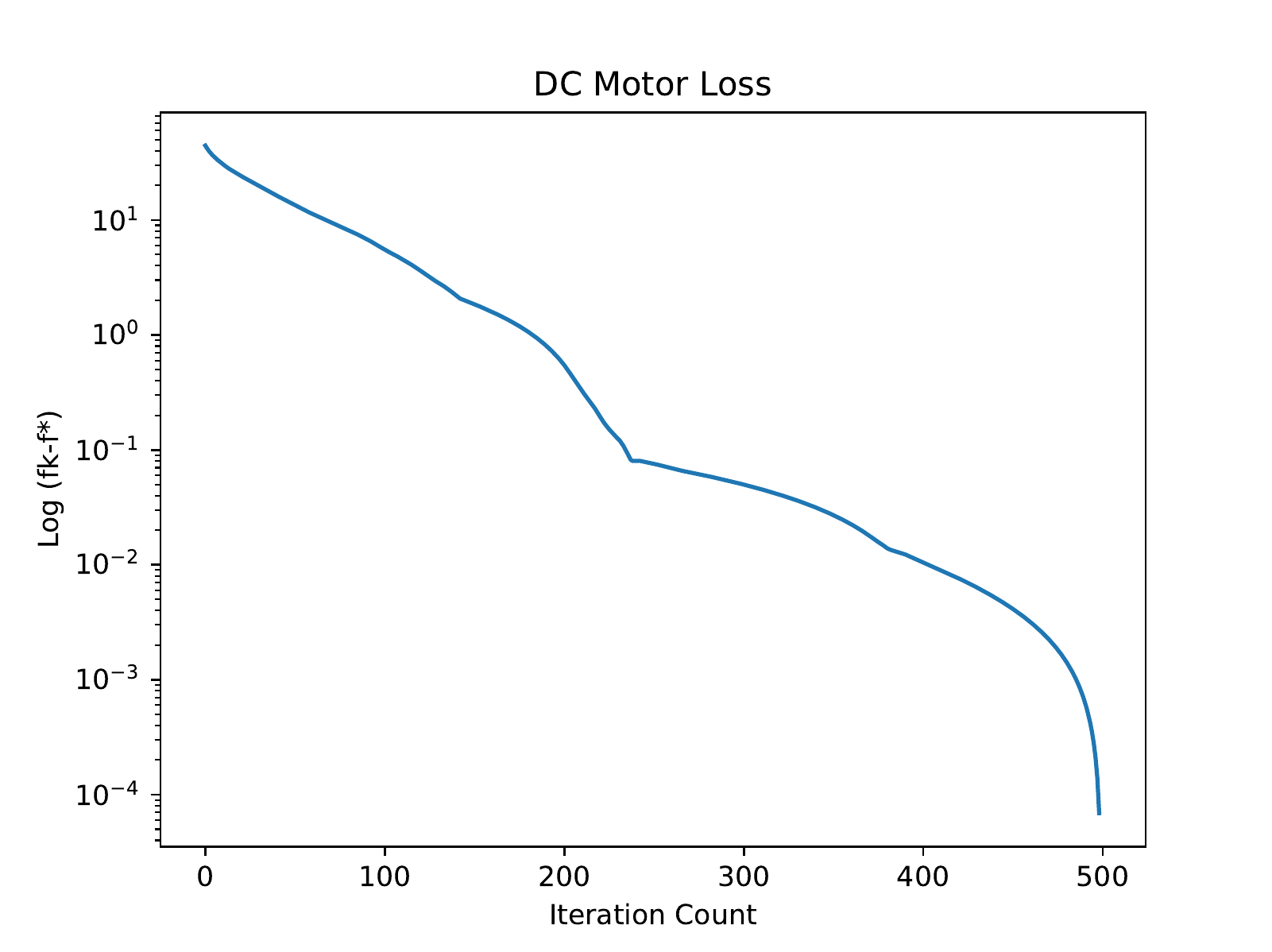} 
    \end{tabular}
 \caption{\label{fig:local_rate} Convergence rate for DRS splitting is locally linear when the objectives are piecewise linear-quadratic (PLQ). The convergence plot show here corresponds to the robust DC motor example in Figure~\ref{fig:DCmotor}.} 
     \end{center}
\end{figure}

\noindent
{\bf Local Linear Rate.} 
When $\rho$ is piecewise linear-quadratic~\cite{RTRW,JMLR:v14:aravkin13a}, the DRS algorithm converges locally linearly to a solution, see Figure~\ref{fig:local_rate}. More precisely, there is a real number $R >0$ such that if $||\eta^K - \eta^*|| < R$ then there is a constant $\kappa$ with $0 < \kappa < 1$ such that for all $k > K$, 
\[
\|\eta^{k+1} - \eta^*\| < \kappa\|\eta^k - \eta^*\|,
\]
where $\eta=\begin{bmatrix}z & \zeta \end{bmatrix}^T$, is the primal and dual pair. 

\begin{theorem}
Algorithm 1 converges with a locally linear rate.
\end{theorem}
\vspace{-.1in}
\textit{ Proof:}
Following the proof technique of ~\cite[Theorem 5]{LFP}, Algorithm 1 has a local linear convergence rate if the following two conditions are satisfied:
\begin{enumerate}
\item Algorithm~\ref{DRS} can be written as the action of a nonlinear operator satisfying a regularity property (see Lemma~\ref{lemma:property}).
\item The functions $g, \rho$ are \textit{subdiffererentially metrically subregular}\footnote{ A mapping $F: \mathbb{R}^n \rightrightarrows \mathbb{R}^m$ is called \textit{metrically subregular} at $\bar{x}$ for $\bar{y}$ if $(\bar{x}, \bar{y}) \in $ graph $F$ and there exists $\eta \in [0 , \infty)$, neighborhoods $\mathcal{U}$ of $\bar{x}$, and $\mathcal{Y}$ of $\bar{y}$ such that 
\[
d(x, F^{-1}\bar{y}) \leq \eta d(\bar{y}, Fx \cap \mathcal{Y}), \quad \forall x \in \mathcal{U}
\]}.
\end{enumerate}
We show that these conditions hold for Algorithm~\ref{DRS}. Define 
\[
Dx \mapsto \begin{bmatrix}\partial g(z) \\ \partial \rho^*(\zeta)\end{bmatrix}, \quad
M = \begin{bmatrix} 0 & I\\ -I & 0\end{bmatrix}, \quad
H = \begin{bmatrix} \frac{1}{\tau}I & 0 \\ -2I & \frac{1}{\sigma}I\end{bmatrix}.
\]

Define the nonlinear operator $T$ by 
\begin{equation}
\label{eq:T}
T = (H+D)^{-1}(H-M).
\end{equation}
$T$ captures the iteration in Algorithm~\ref{DRS}, which can be written as $\eta^{k} = T\eta^{k-1}$, for 
$\eta = \begin{bmatrix} z^T, \zeta ^T\end{bmatrix}^T$. Then we have the following lemma.
\begin{lemma}
\label{lemma:property}
Suppose that $\tau, \sigma < 1$. Then 
\[||T\eta-\eta||_{H-M}^2 \leq \langle \eta^*-\eta,(H-M)(T\eta-\eta)\rangle
\]
where $\eta^*$ is such that $0 \in (D+M)\eta^*$.
\end{lemma}
The proof is given in the Appendix. \\
\begin{figure}
  \begin{center}
 \includegraphics[scale=0.6]{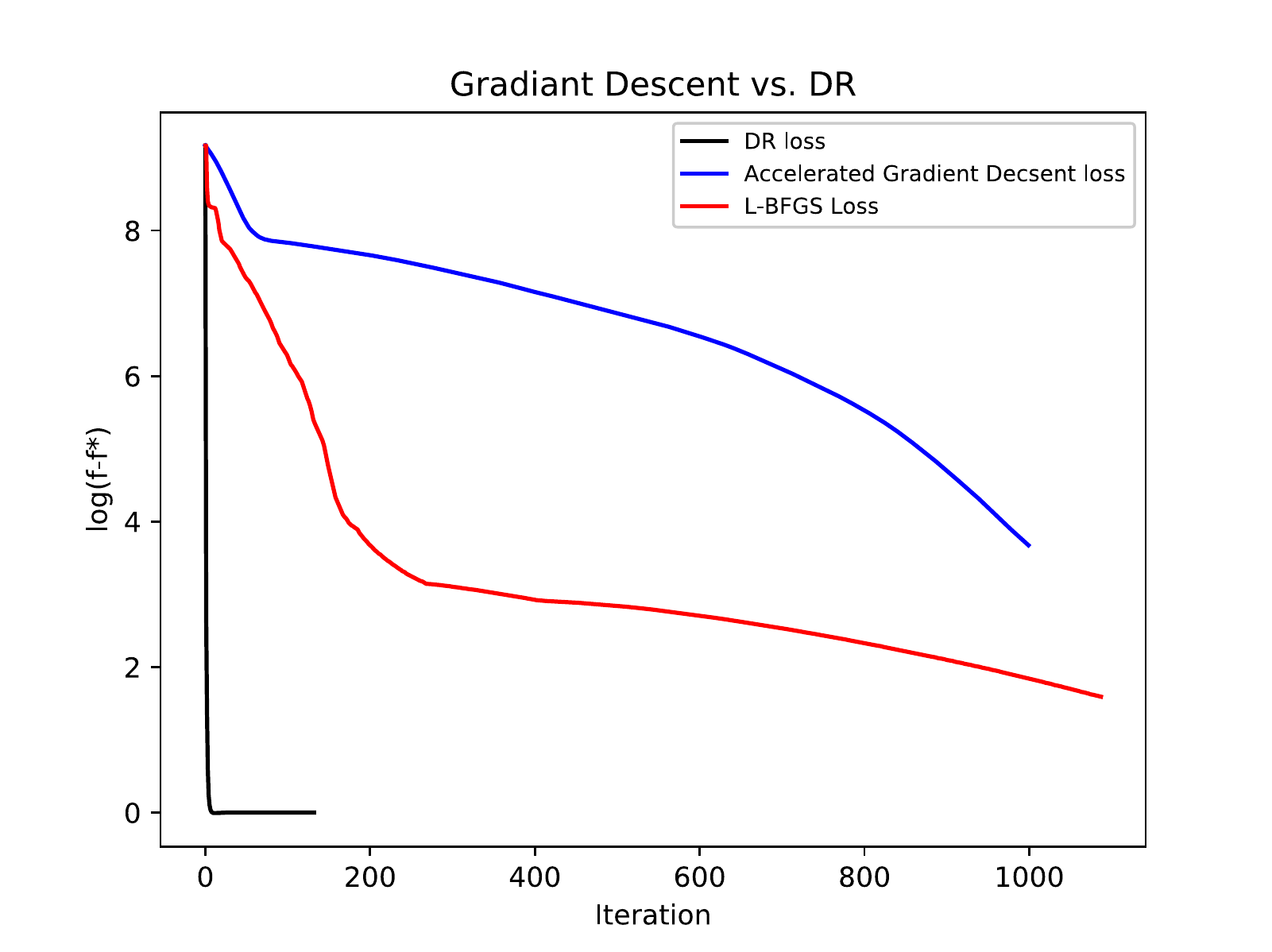} 
 \caption{\label{fig:comparison} Objective vs. iteration counts of Algorithm~\ref{DRS} for~\eqref{eq:full} (black),
 vs. accelerated gradient descent (AGD) (blue) and L-BFGS (red) for~\eqref{eq:genSmoother}. 
 Both $\rho_1$ and $\rho_2$ are Huber losses, with $\rho_3 \equiv 0$, $N=200, n=2$ and $Q, R$ nonsingular, so~\eqref{eq:full} and~\eqref{eq:genSmoother}  are equivalent. All iterations require $O(n^2N)$ operations. 
 DRS splitting is much faster than methods with linear convergence rates and similar iteration complexity. } 
     \end{center}
\end{figure}
\begin{figure}
  \begin{center}
 \includegraphics[scale=0.4]{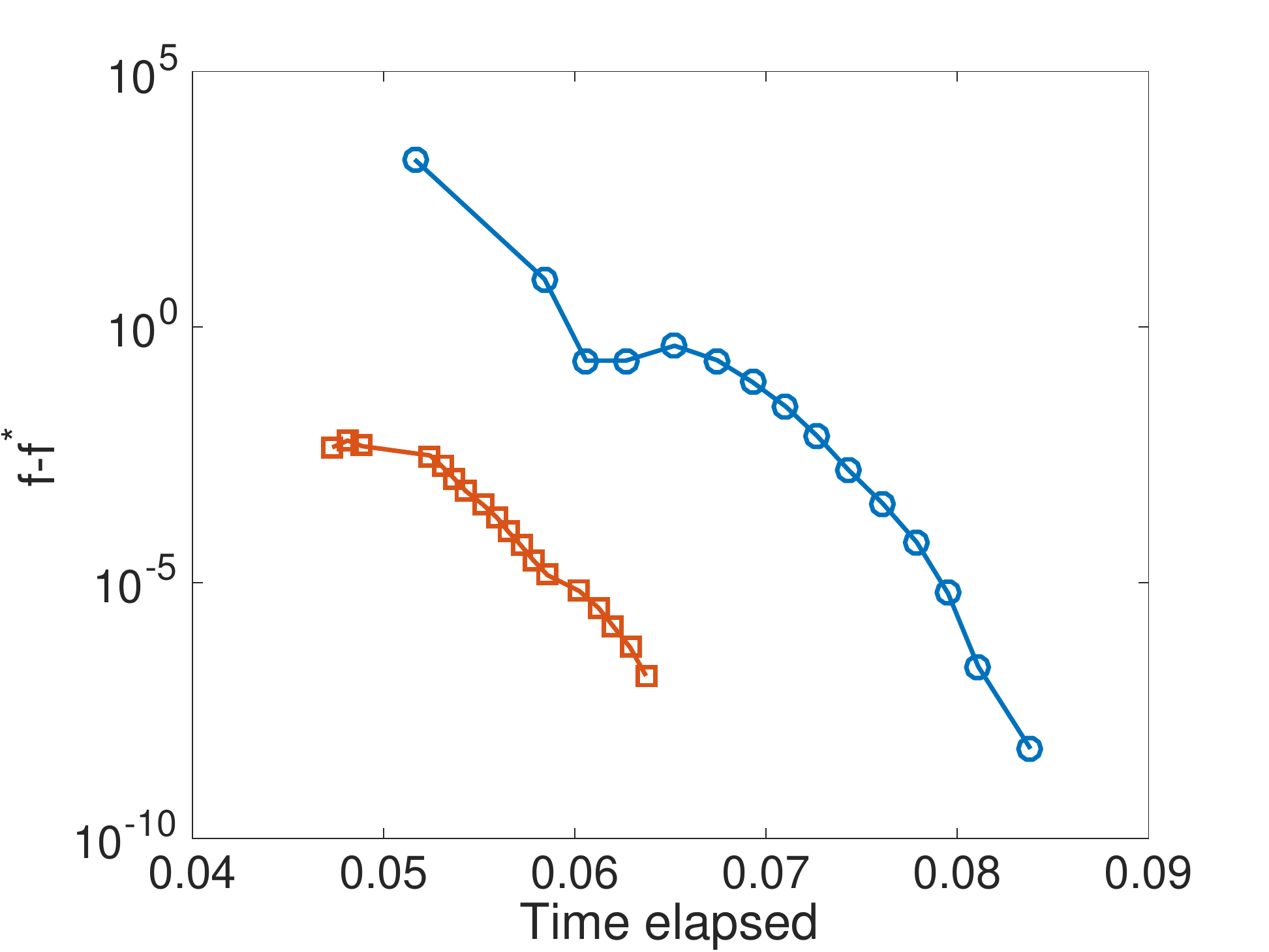}
 \caption{\label{fig:timecompare} Timed run of Algorithm~\ref{DRS} vs. IPsolve for the same setup as presented in Figure~\ref{fig:comparison}. 
At this scale, we see the locally linear convergence rate of the DRS. Even though IPsolve has superlinear rate, DRS wins because the slope 
of the rate is very steep, and each iteration is fast. By the time DRS is done, IPsolve has had time for only taken a few iterations.}
  \end{center}
\end{figure}
This establishes condition (1). 
Condition (2) requires the concept of metric subregularity. This property holds for PLQ functions~\cite{LFP}, 
and holds for indicators of convex sets by~\cite[Theorem 3.3]{AG},  reproduced in the Appendix. This completes the proof of the theorem.

\noindent
{\bf Comparison on Smooth Nonsingular Problems.}
If the covariances, $Q, R$ are non-singular and the penalties $\rho_{1,2}$ are $\mathcal{C}^1$-smooth, 
then the Kalman smoothing problem can be written as a smooth convex problem. In this case the same reformulation will work and Algorithm~\ref{DRS} will still give a local linear rate. However more common algorithms such as gradient descent and L-BFGS can also be applied. We compare the performance of these three algorithms to track a particle moving along a smooth path with $N=200$ and $n=2$.
We use non-singular $Q_k$, and Huber penalty functions.\\
As seen in Figure~\ref{fig:comparison}, Algorithm~\ref{DRS} for~\eqref{eq:full} converges far faster than either accelerated gradient descent 
or LBFGS method on the equivalent nonsingular smoothing formulation~\eqref{eq:genSmoother} . This is because its convergence rate does not 
depend on the condition number of the matrix $A$, so each iteration makes a lot of progress, and 
we can keep the complexity of each iteration at $O(n^2N)$, same as for a matrix-vector multiply needed for a gradient evaluation, 
if we factor the sparse block tridiagonal matrix $AA^T$ once at the start of the algorithm.  \\
We also compare with the second-order interior point method, implemented in the IPsolve package\footnote{\url{https://github.com/UW-AMO/IPsolve}.}. Use-cases and performance of IPsolve for nonsingular Kalman smoothing is discussed in~\cite{aravkin2017generalized}. The results are shown in Figure~\ref{fig:timecompare}, where IPsolve and 
DRS for the equivalent reformulation are compared for the nonsingular Huber model.
Even though DRS has at best a linear rate, the constants are very good, as they do not depend on the conditioning 
of the Kalman smoothing problem.
 The other advantage is that DRS can use a pre-factorized matrix, while IPsolve has to solve a modified linear system 
 every time; there is no simple strategy to pre-factor as with DRS.\\
The numerical experiments suggest that Algorithm~\ref{DRS} should be used regardless of whether $Q$ and $R$ are singular or not.  In the next section, we focus on a rich class of singular noise models found in navigation.

%% file: vehicleModel_paper.tex
\section{Navigation Models}
\label{sec:vehicleModel}

\vspace{-.1in}

\begin{figure}[h!]
\begin{center}
\includegraphics[scale=1]{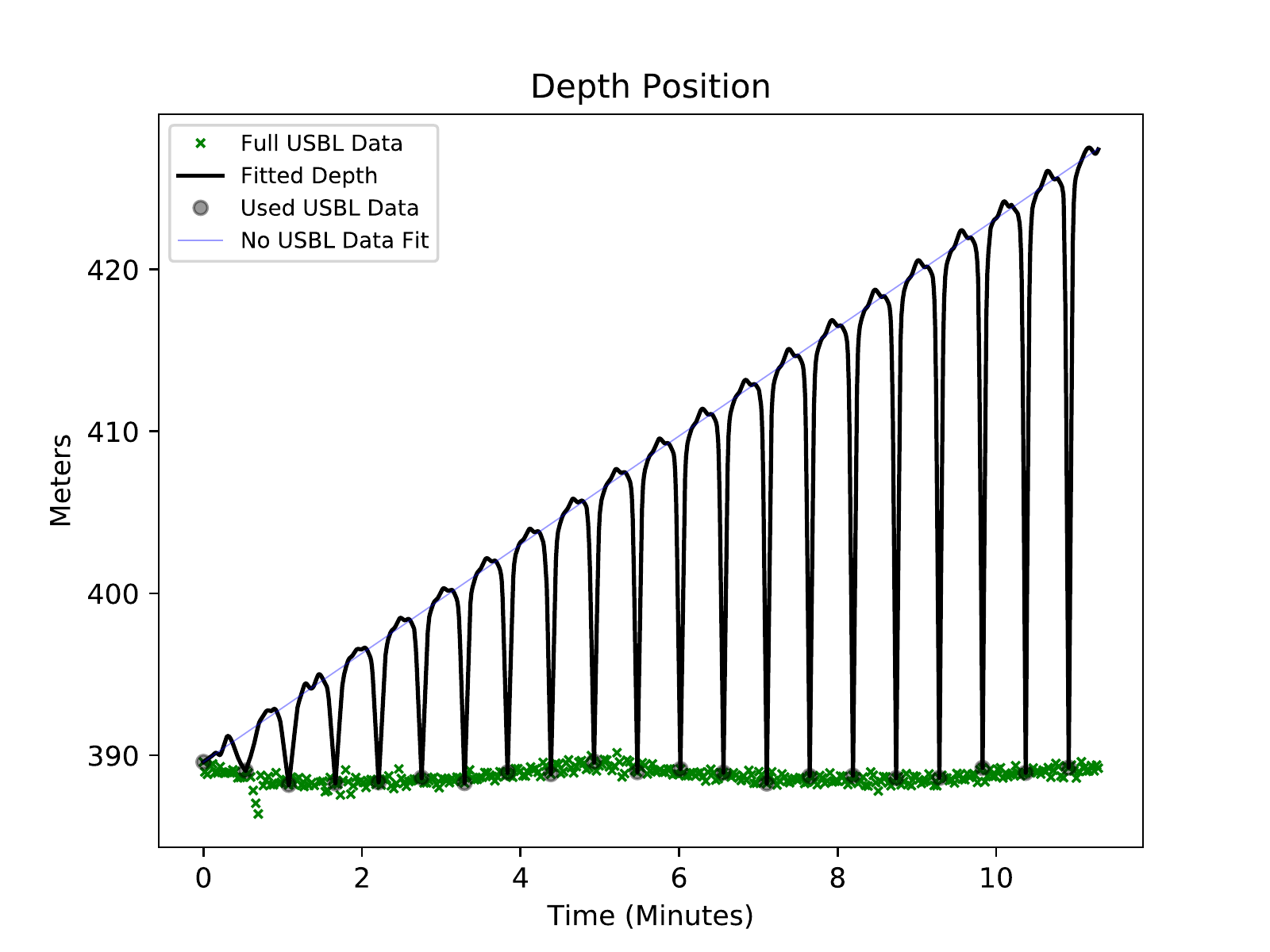}
\caption{Full position data is plotted using green crosses; the smoother only uses 
a subset highlighted with gray disks. The model is used to provide the
initialization for the solver, shown in blue, by forward propagating
from the first data point. With $Q$ as in~\eqref{eq:Qderiv}, the model forces a zero
acceleration constraint. Combined with a non-zero initial velocity, this results
in an erroneous initialization, which even high
confidence in the observed datapoints is unable to overcome, yielding a counterintuitive 
result.}
\label{fig:badexample}
\end{center}
\end{figure}

Autonomous navigation requires high-fidelity tracking using occasional GPS and frequent depth/height, gyrocompass, and linear acceleration data. Gyro, compass, and linear acceleration are readily available from inertial measurement units (IMUs). \\
In this section, we develop a simple kinematic model that is trivially applicable to any
vehicle, and is particularly  appropriate for many underwater vehicle applications, where
accelerations are heavily damped and autonomous vehicles often travel in long straight lines (e.g. for survey work). When the attitude is known or changing slowly, the model can be linearized effectively and the situation simplifies considerably; 
our synthetic examples and underwater survey application use linearized models. \\
{\bf Linear Singular Navigation Model.}
For a vehicle that is well-instrumented in attitude, the uncertainty in position
(and the x-y states in particular) is typically orders of magnitude larger than
the uncertainty in attitude. In practice, we simplify the full
nonlinear vehicle process model to track only position states 
$( x, y, z)$, while assuming that
the attitude states $(r, p,h)$ are directly available from the most recent sensor
measurements. To make the model linear, the position and its derivatives are
referenced to the local-level frame.\\
To incorporate linear acceleration measurements from
an inertial measurement unit (IMU), we must track both linear velocities 
and linear acceleration in the state vector.  This leads to the augmented state
\begin{equation}
  \vect{x}_s = [x,y,z,\dot{x},\dot{y},\dot{z},\ddot{x},\ddot{y},\ddot{z}]^\top.
\end{equation} \label{lin_pm}
The linear kinematic process model is given by 
\begin{align}
\dot{\vect{x}}_s &= 
\underbrace{\left[\begin{array}{ccc}
\mtrx{0} & \mtrx{I} & \mtrx{0} \\
\mtrx{0} & \mtrx{0} & \mtrx{I} \\
\mtrx{0} & \mtrx{0} & \mtrx{0} \end{array}\right]}_{\mbox{$\mtrx{F}_s$}} \vect{x}_s + 
\underbrace{\left[\begin{array}{c}
\mtrx{0} \\
\mtrx{I} \\
\mtrx{0} \end{array}\right]}_{\mbox{$\mtrx{G}_s$}} \vect{w}_s \label{x_s-dot},
\end{align}
where $\vect{w}_s \sim \mathcal{N}(0,\mtrx{Q}_s)$ is zero-mean
Gaussian  noise. \\
The linear process model \eqref{x_s-dot} is usually 
discretized using a Taylor series: 
\begin{align}
\vect{x}_{s_{k+1}} &= \mtrx{F}_{s_k} \vect{x}_{s_k} + \vect{w}_{s_k} \label{disc_s}\\
\mtrx{F}_{s_k} &= e^{\mtrx{F}_s T}\\
 &= \mtrx{I} + \mtrx{F}_sT +{\frac{1}{2!}\mtrx{F}_s^2T^2}  + \cancelto{0}{\frac{1}{3!}\mtrx{F}_s^3T^3} + \cdots \nonumber\\
 &= \left[\begin{array}{ccc}
\mtrx{I}&\mtrx{I}T&\frac{1}{2}\mtrx{I}T^2\\
\mtrx{0}&\mtrx{I}&\mtrx{I}T\\
\mtrx{0}&\mtrx{0}&\mtrx{I} \end{array} \right] \nonumber
\end{align}
where the higher order terms are {\it identically zero} because of the
structure of $\mtrx{F}_s$, resulting in a simple closed-form solution
for $\mtrx{F}_{s_k}$. 
%
The discretized process noise
\begin{equation}
  \vect{w}_{s_k} = \int^T_0 e^{\mtrx{F}_s(T - \tau)}\mtrx{G}_s \vect{w}_s(\tau) d\tau,
\end{equation}
is a zero-mean Gaussian, with covariance given by 
\begin{equation}
\label{eq:Qderiv_0}
\mtrx{Q}_{s_k} = \int^T_0 e^{\mtrx{F} (T-\tau)} \mtrx{G} \mtrx{Q} \mtrx{G}^\top
e^{\mtrx{F}^\top (T-\tau)} d\tau, 
\end{equation}

which simplifies to
\begin{equation}
\label{eq:Qderiv}
\mtrx{Q}_{s_k} = 
\left[\begin{array}{ccc}
\frac{1}{3}T^3&\frac{1}{2}T^2&0 \\
\frac{1}{2}T^2&T             &0 \\
0             &0             &0 \end{array} \right] \mtrx{Q}_s,
\end{equation}
for
\begin{equation*}
e^{\mtrx{F} (T-\tau)} = \left[\begin{array}{ccc}
\mtrx{I}&\mtrx{I}(T-\tau)&\frac{1}{2}\mtrx{I}(T-\tau)^2\\
\mtrx{0}&\mtrx{I}&\mtrx{I}(T-\tau)\\
\mtrx{0}&\mtrx{0}&\mtrx{I} \end{array} \right],\ \ \mtrx{G} = \left[\begin{array}{c}
\mtrx{0}\\
\mtrx{I}\\
\mtrx{0} \end{array} \right].
\end{equation*}
In practice this can lead to wildly incorrect results. In Figure~\ref{fig:badexample}, 
we show the estimate of position obtained from a subset of the navigation data.
The model is defined with constraints
\[
Q_{s_k}^{1/2}u_{s_k} = F_{s_k}x_{s_{k+1}} - x_{s_k}.
\]
The $Q$ in~\eqref{eq:Qderiv} forces the acceleration to be $0$ across the entire model because the lower right corner 
is set to $0$.  As a result, the initialized track can be biased away from the
data by a fixed
velocity, obtained by finding the slope from the most recent position data. The available data do not agree, but the constraint is stronger; the
information is integrated in a counter-intuitive way. \\
Instead, we model the covariance as if the error were the next term in the Taylor series approximation, 
a technique suggested by~\cite{YAA}. More precisely we set covariance to be the outer product, $\Gamma^T \Gamma$ where
\[
\Gamma = \begin{bmatrix}\frac{1}{3!}\mtrx{I}T^3 & \frac{1}{2!}\mtrx{I}T^2 & \mtrx{I}T\end{bmatrix}
\]
This leads to a rank 3 covariance for a $9\times 9$ matrix for a model that comprises
position, velocity, and acceleration in 3D space. This model avoids the issue in Figure~\ref{fig:badexample}.\\
{\bf Measurement Models for the IMU.}
The inertial measurement unit (IMU) does not measure position or velocity, just linear and angular accelerations. 
To use these measurements, we track linear acceleration as part of the state. However, the 
acceleration measured by the IMU is relative to the physical frame of the
vehicle on which it is mounted, while the acceleration of the state 
is relative to the navigation frame. A coordinate transformation between these frames is required for a comparison; 
we use heading, pitch, and roll of the vehicle for the linear
model.
The transformation from body-frame to local-level is given by \textbf{$R(\vect{\varphi})$},
where $\vect{\varphi}$ comprises heading $h$, pitch $p$, and roll $r$:
\begin{equation}
  \mtrx{R}(\vect{\varphi}) = \mtrx{R}^\top_h \mtrx{R}^\top_p \mtrx{R}^\top_r,
\end{equation}
where $R_h$, $R_p$, and $R_r$ are given by 
\begin{equation}
\left[\begin{array}{ccc}
      c{h} & s{h} & 0 \\
      -s{h} & c{h} & 0 \\
      0 & 0 & 1 \end{array}\right],\quad  
      \left[\begin{array}{ccc}
      c{p} & 0 & -s{p} \\
      0 & 1 & 0 \\
      sp & 0 & cp \end{array}\right], \quad 
\left[\begin{array}{ccc}
      1 & 0 & 0 \\
      0 &c{r} & s{r} \\
      0 & -s{r} & c{r} \end{array}\right]
\end{equation}
with $c\cdot$ and $s\cdot$ shorthand for $\cos(\cdot)$ 
and $\sin(\cdot)$.  \\
Any navigation system that relies on an IMU needs occasional measurements that inform the position 
(e.g. GPS), otherwise the error in position estimates grows without bound. We are given 
these data from a separate source, sampled at a lower update rate than that of the IMU. 
For any $s$ where such data is available, we have the measurement model 
\[
H_s = \begin{bmatrix} I_{3 \times 3} & 0_{3 \times 6}\\ 0_{3 \times 6} & R(\vect{\varphi})\end{bmatrix}, \quad
\vect{z}_s = \begin{bmatrix} \vect{b}^\top & \ddot{x}_{\textrm{meas}} &  \ddot{y}_{\textrm{meas}} & \ddot{z}_{\textrm{meas}}\end{bmatrix}^\top.
\]
If there is no position data measured at time $s$ then we use the model 
\[
H_s = \begin{bmatrix} 0_{3 \times 3} & 0_{3 \times 6}\\ 0_{3 \times 6} & R(\vect{\varphi})\end{bmatrix}, \quad 
\vect{z}_s = \begin{bmatrix} 0& \ddot{x}_{\textrm{meas}} & \ddot{y}_{\textrm{meas}} & \ddot{z}_{\textrm{meas}}\end{bmatrix}^\top.
\]
The covariance used for measurement data depends on whether there was position data available:
\[
R_s = \begin{bmatrix}  0_{3 \times 3} & 0_{3 \times 3}\\ 0_{3 \times 3} & r_s I_{3 \times 3}\end{bmatrix}, \quad R_s = \begin{bmatrix} U_s & 0_{3 \times 3}\\ 0_{3 \times 3} & r_s I_{3 \times 3}\end{bmatrix}
\]
wher ethe top $3\times 3$ block is either $0$ (position not available) or $U$, a diagonal matrix reflecting position uncertainty (position is available). 
The scalar $r_s$ models uncertainty in IMU measurements. 



%% file: experiment.tex
\section{Analysis of Mooring Data}
\label{sec:numerics}

We are interested in the ability to maintain an accurate position
estimate on-board an autonomous underwater vehicle using acceleration
measurements from a low cost  inertial measurement unit (IMU), given periodic
position fixes. To test this, we use the singular general Kalman framework to analyze data collected from a surface mooring
equipped with an IMU that was deployed off the coast of Florida during spring 2017. 
We use the mooring, which is drifting with the current, as a proxy for a slowly
moving underwater vehicle subject to unknown disturbances. 
In particular we are looking at the position uncertainty and error accrued over time between the
periodic, world-referenced position fixes that are provided by the ultra short
baseline (USBL) system.

The new capabilities are useful because
\begin{enumerate}
\item Navigation models are singular
\item Data are noisy
\item IMU has biases, captured using singular models
\item data can be quantized, motivating a special loss.
\end{enumerate}
 In this analysis, we use the singular linear kinematics model
in Section \ref{sec:vehicleModel}, the Huber loss from Section~\ref{sec:intro},
and the DRS algorithm from Section~\ref{sec:Opt} to solve the final smoothing
problem.

\subsection{Experimental Setup}
\label{sec:exp_setup}

\vspace{-.1in}

\begin{figure}[h!]
  \begin{center}
    \includegraphics[width=\columnwidth]{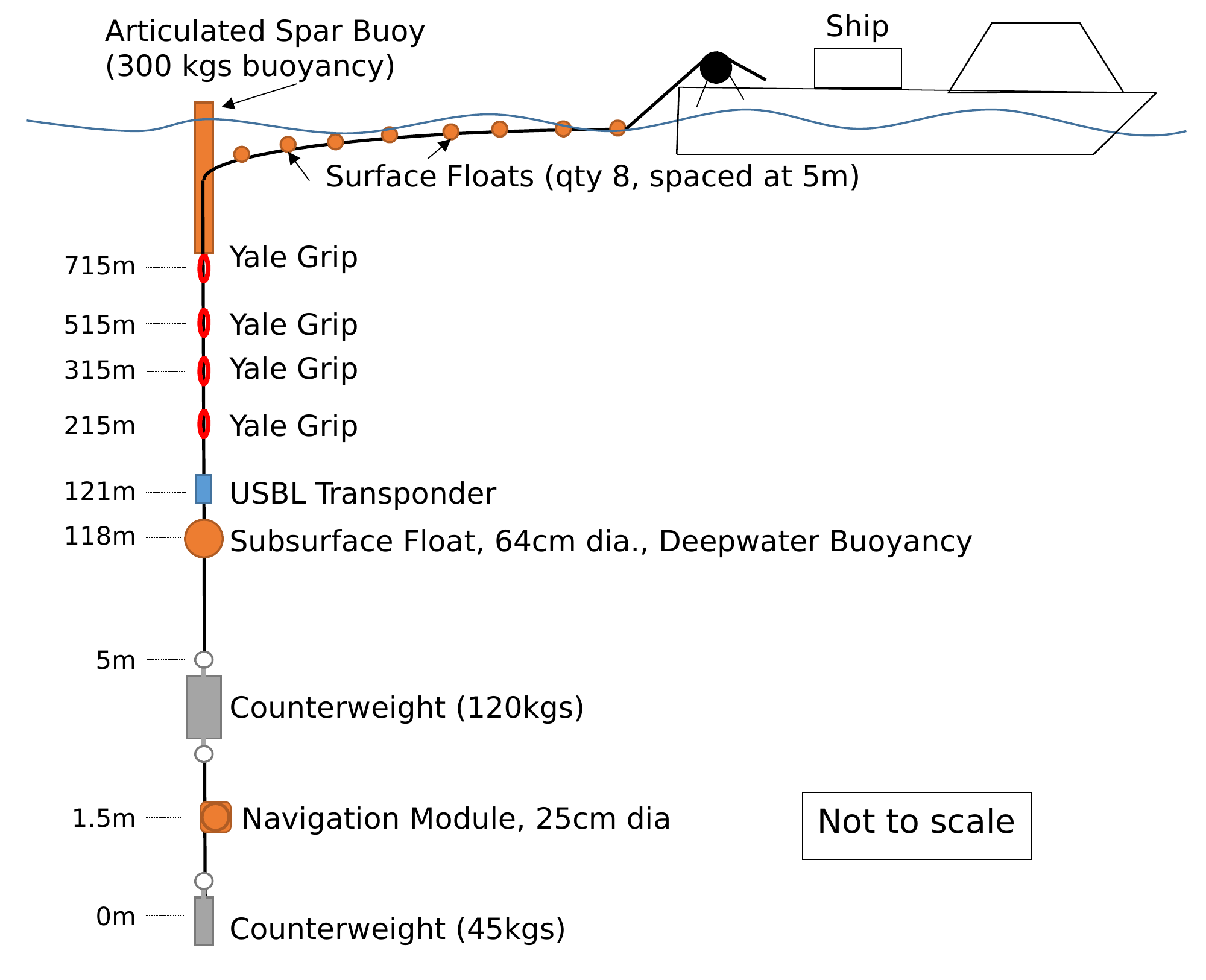}
    \caption{Sketch of the mooring and tender ship showing major buoyancy and
      counterweight components and their relative depths. The USBL system
      provides position fixes of the USBL transponder. The
      navigation module provides roll, pitch, heading, and linear accelerations.}
  \label{fig:depl_setup}
  \end{center}
\end{figure}

As shown in Figure \ref{fig:depl_setup}, the
mooring comprises  an articulated spar buoy on the surface, supporting a cable
with various instruments attached. The mooring can be shortened using
yale grips shown in figure. We are using a portion of data from when the mooring
is at its max length of approximately 715m.
At 121 meters above the bottom of the mooring is an ultra short baseline (USBL)
receiver which, in concert with a nearby tender ship, provided three dimensional
position updates for the mooring (latitude, longitude, and depth). 
Below the main clump weight is a 4.25 meter section of
Spectra\textsuperscript{\textregistered} line with its own smaller clump weight
of approximately 45 kg. This supports a 25 cm diameter spherical glass housing
containing the navigation module.\\
\begin{figure}[h]
  \begin{center}
    \includegraphics[width=0.7\columnwidth]{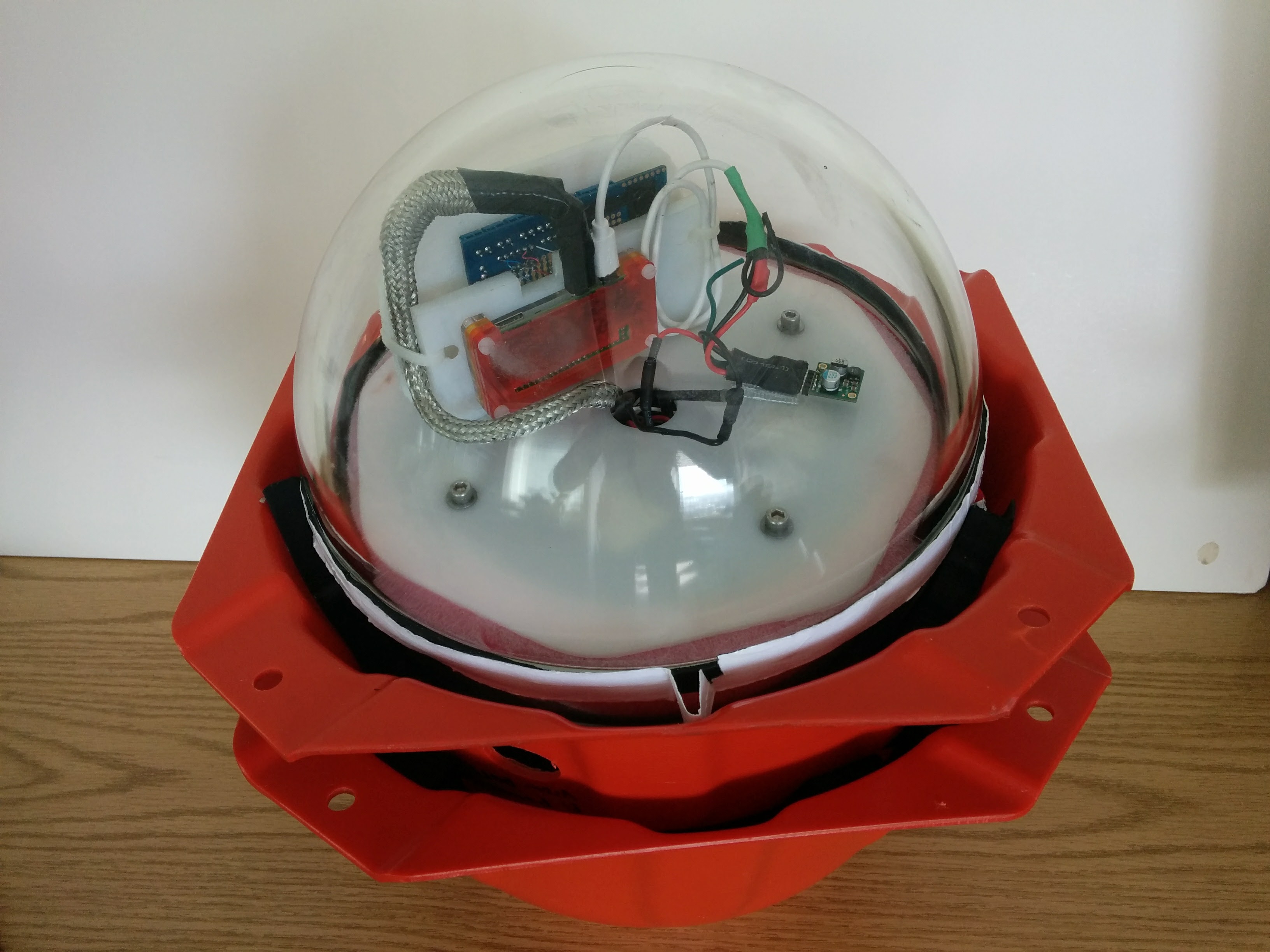}
    \caption{The navigation module, housed in a 25 cm diameter glass sphere. 
      }
  \label{fig:navmodule}
  \end{center}
\end{figure}
The self-contained navigation module, shown in Figure \ref{fig:navmodule},
consists of a RaspberryPi-based logger supporting a precision clock (Adafruit
ChronoDot RTC v2.1, based on the DS3231 temperature compensated crystal
oscillator), gyro (L3GD20H), and accelerometer and compass (LSM303D). The
navigation module carries its own batteries and recorded continuously throughout the
deployment, providing the time-stamped attitude and acceleration data used in this analysis.
Quantization in the attitude (roll, pitch, and yaw) and
linear acceleration measurements resulted in a degradation of the native
accuracy of the sensors. Table  \ref{tbl:specs} provides a summary of the
measurements and associated resolutions as recorded during this experiment.
In this capacity the navigation module data serves as a proxy for a low cost autonomous
underwater vehicle using a low grade commercial IMU.\\
\begin{table}[ht]
\renewcommand{\arraystretch}{0.97} \scriptsize
  \begin{tabular}{p{0.25\columnwidth}p{0.22\columnwidth}p{0.25\columnwidth}}
  \hline
  \textbf{Measurement} & \textbf{Resolution} & \textbf{Sample Freq.} \\
  \hline
  time                           & 3.5 ppm   & n/a\\
  roll, pitch, yaw               & 0.1$^{\circ}$   & 25 Hz\\
  lin. acceleration              & 0.00766 m/s$^2$   & 25 Hz\\
  \hline
  \end{tabular}
\caption{~Navigation module sensor specifications.}
\label{tbl:specs}
\end{table}

The articulated spar buoy was tethered to the ship through an umbilical that
supplied power, two way communications and data transfer. During operations,
the intent was to decouple the motion of the ship from that of the surface
mooring, keeping slack in this umbilical. This is accomplished by using the ship
to tow the mooring into position and then allowing both the ship and mooring to
drift with the current. 

Ground truth for the position of the mooring was provided by the Sonardyne
Ranger 2, a USBL system that provided 3-D position fixes every 2 seconds. The
USBL system self reports its measurement uncertainties at each
measurement. These ranged from 3.7 to 7.5 m uncertainty in x and y, and 0.8 to 4.0 m uncertainty in depth.

\subsection{Model and Experimental Results}

\vspace{-.1in}




%

\begin{figure}[h!]
\begin{center}
\includegraphics[scale=.5]{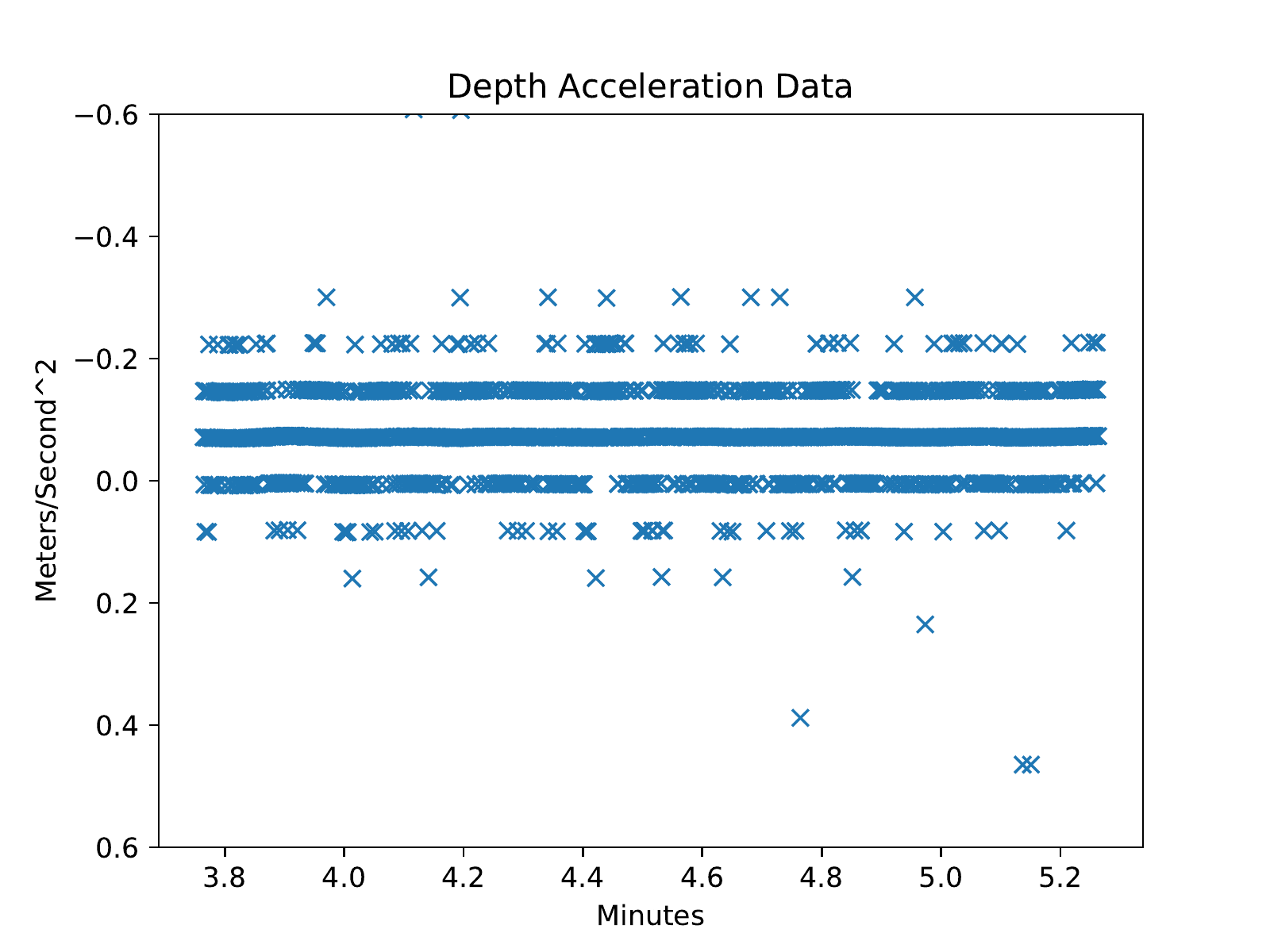}
\caption{A snippet of the depth acceleration data, rotated into the world frame,
  shows the quantization of the acceleration data.}
\label{fig:acceldata}
\end{center}
\end{figure}

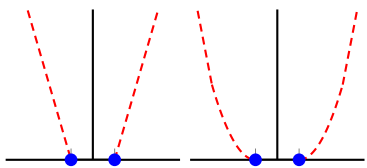
\begin{figure}[h!]
\begin{center}
\begin{tikzpicture}
  \begin{axis}[
    thick,
    height=2cm,
    xmin=-2,xmax=2,ymin=0,ymax=1,
    no markers,
    samples=50,
    axis lines*=left, 
    axis lines*=middle, 
    scale only axis,
    xtick={-0.5,0.5},
    xticklabels={},
    ytick={0},
    ] 
    \addplot[red,domain=-2:-0.5,densely dashed] {-x-0.5};
    \addplot[domain=-0.5:+0.5] {0};
    \addplot[red,domain=+0.5:+2,densely dashed] {x-0.5};
    \addplot[blue,mark=*,only marks] coordinates {(-0.5,0) (0.5,0)};
  \end{axis}
\end{tikzpicture}
\begin{tikzpicture}
  \begin{axis}[
    thick,
    height=2cm,
    xmin=-2,xmax=2,ymin=0,ymax=1,
    no markers,
    samples=50,
    axis lines*=left, 
    axis lines*=middle, 
    scale only axis,
    xtick={-0.5,0.5},
    xticklabels={},
    ytick={0},
    ] 
    \addplot[red,domain=-2:-1.5,densely dashed] {-1.5*x-1.75};
    \addplot[red,domain=-1.5:-.5,densely dashed] {0.5*(x+0.5)^2};
    \addplot[domain=-0.5:+0.5] {0};
    \addplot[red,domain=+0.5:+1.5,densely dashed] {0.5*(x-0.5)^2};
    \addplot[red,domain=1.5:2,densely dashed] {1.5*x-1.75};
    \addplot[blue,mark=*,only marks] coordinates {(-0.5,0) (0.5,0)};
  \end{axis}
\end{tikzpicture}
\end{center}
\caption{Vapnik loss function and a smoothed variant.}
\label{fig:vap}
\end{figure}

\begin{figure}[h!]
\begin{center}
\includegraphics[scale=.5]{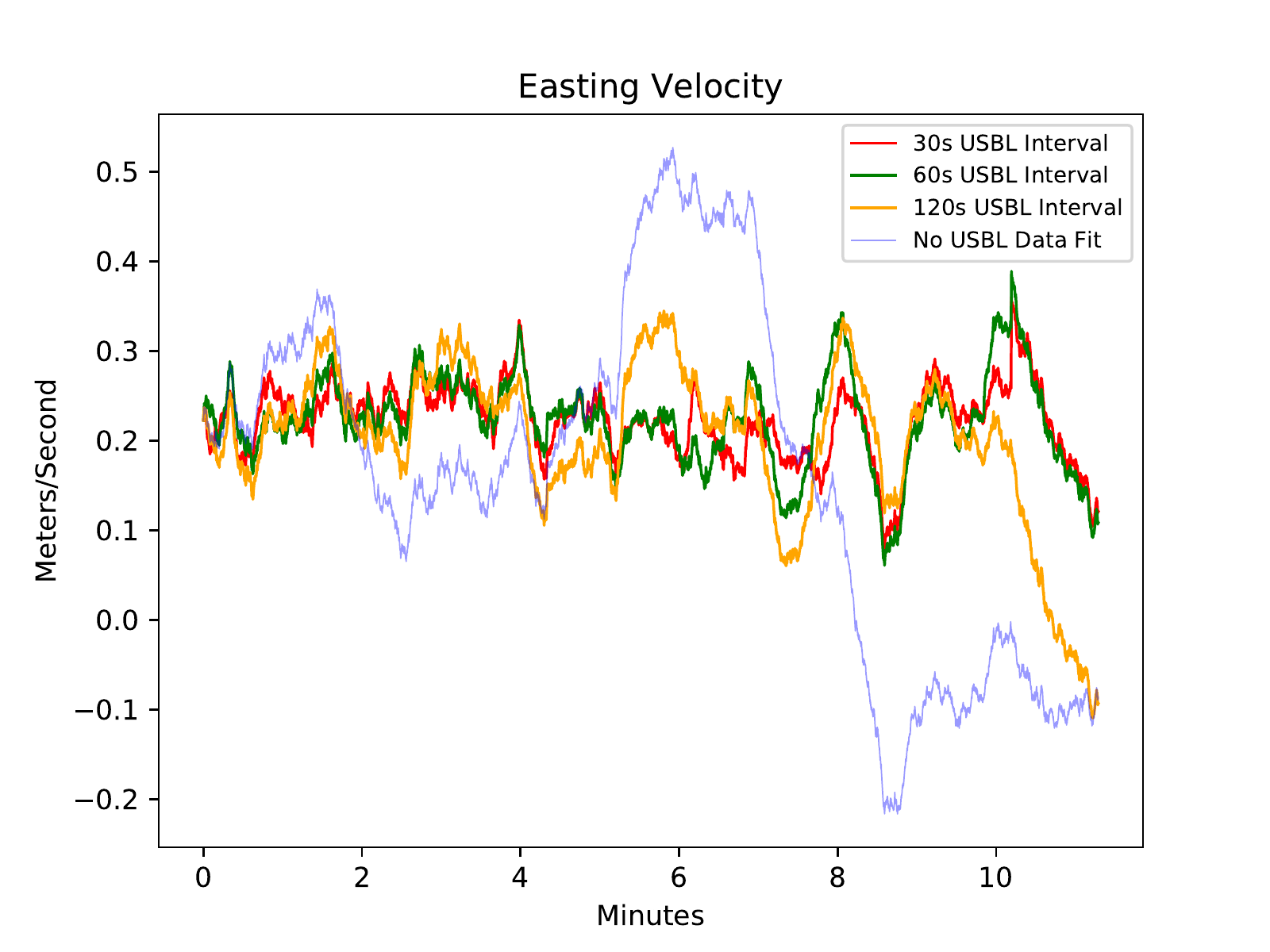}
\includegraphics[scale=.5]{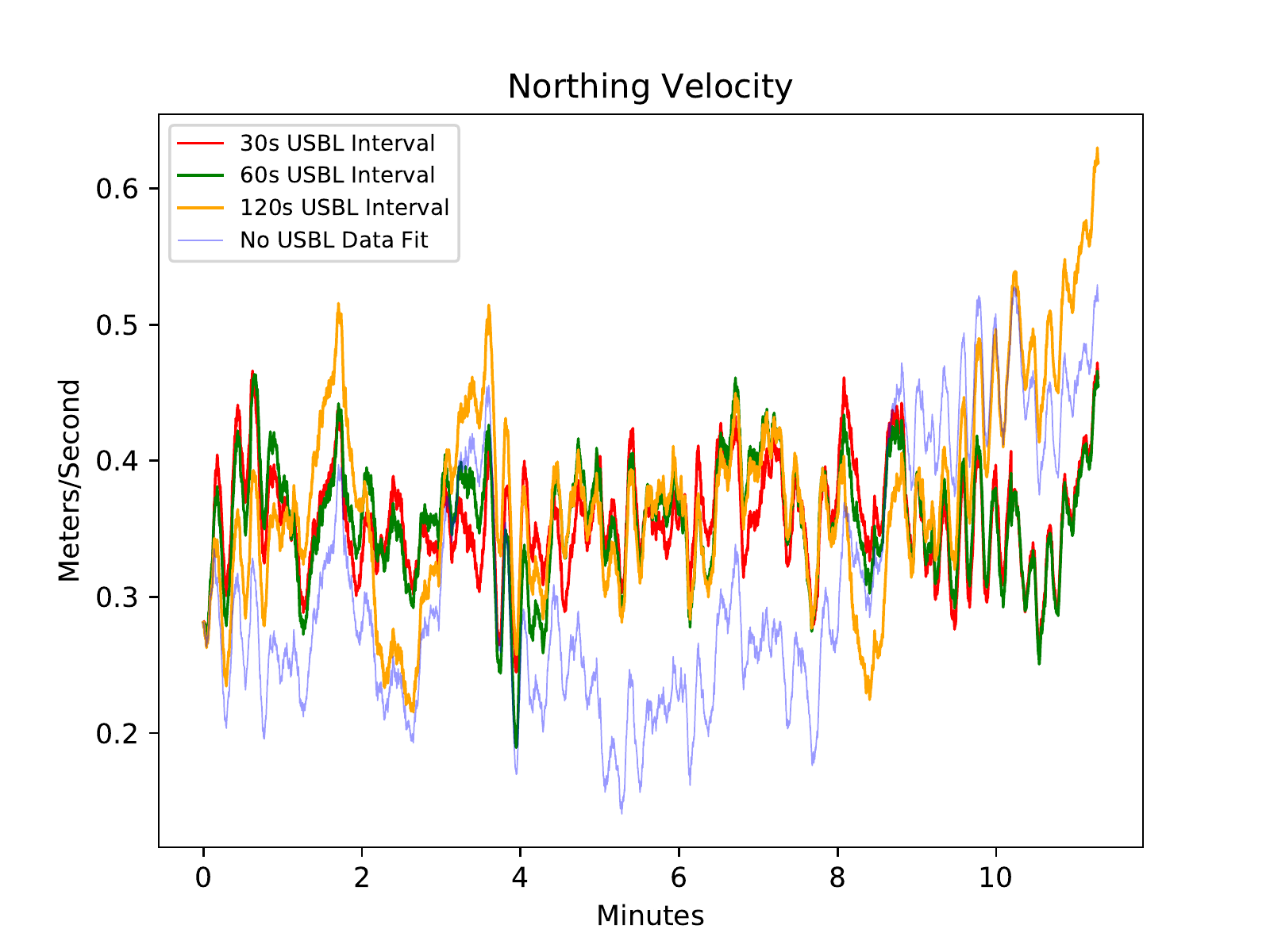}
\includegraphics[scale=.5]{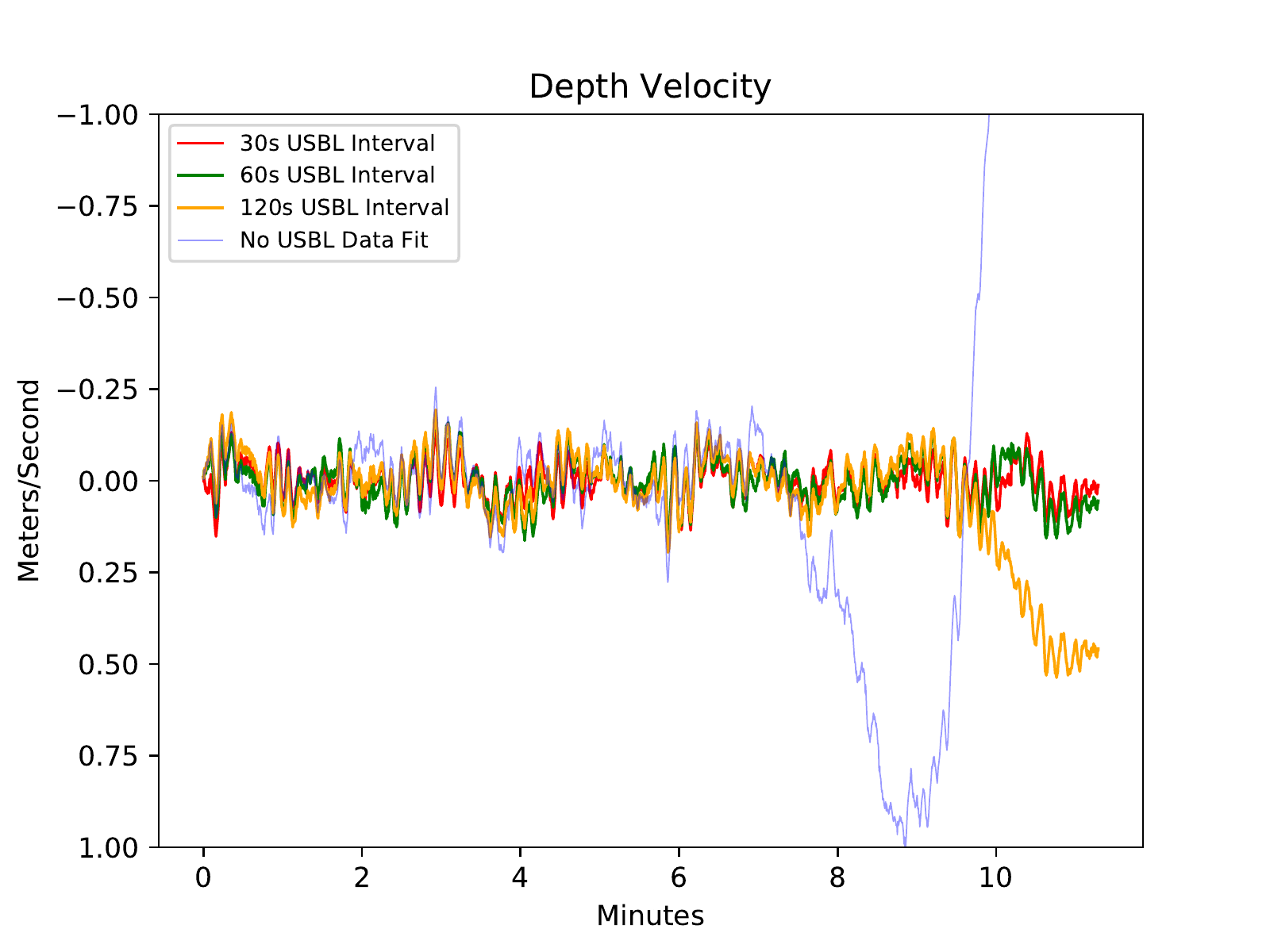}
\caption{Velocity fit for different frequencies of position data. The build up of acceleration errors can be seen clearly here, especially when no position data is used.}
\label{fig:velshort}
\end{center}
\end{figure}

Two challenges in the experimental setup required the flexibility of the modeling framework. 
The depth acceleration data, some of which is plotted in Figure~\ref{fig:acceldata}, 
is extremely discretized and appears to have mean shifted away from zero.
To counteract this, a constant bias for acceleration measurements was fit and removed. 
In the singular framework, we easily include a constant term, by imposing equality constraints 
across all time points using the process model. The measurement maps are then modified 
to directly subtract the estimated bias. \\
Because of the level of discretization we want to 
use the Vapnik loss function (Figure~\ref{fig:vap}) that does not penalize in a small interval around the data. 
The `deadzone' region is set according to the quantization of the data, which is $.05$. 
The `corners' of the Vapnik encourage the errors to be exactly equal to the quantization value, an unnecessary artifact. 
We therefore use a Huberized version of the Vapnik, smoothing the corners but leaving the deadzone. 
In addition to the deadzone, this loss is robust, as it has linear tail growth. \\
{\bf Results for $10$ Minute Track:}
We begin by considering $10$ minutes of IMU data with occasional USBL position data. The position data are available approximately every $2$ seconds, but we test performance with intervals of $30, 60, 120$ seconds. 
The $x_0$ given to the algorithm is as follows: position is set to the first
position fix and acceleration is set to zero, while velocity is taken 
to be the slope from the last available position data to the starting time. 
The algorithm is initialized by propagating this $x_0$ through the entire model and then run for $500$ iterations. \\
\begin{figure}[h!]
\begin{center}
\includegraphics[scale=.5]{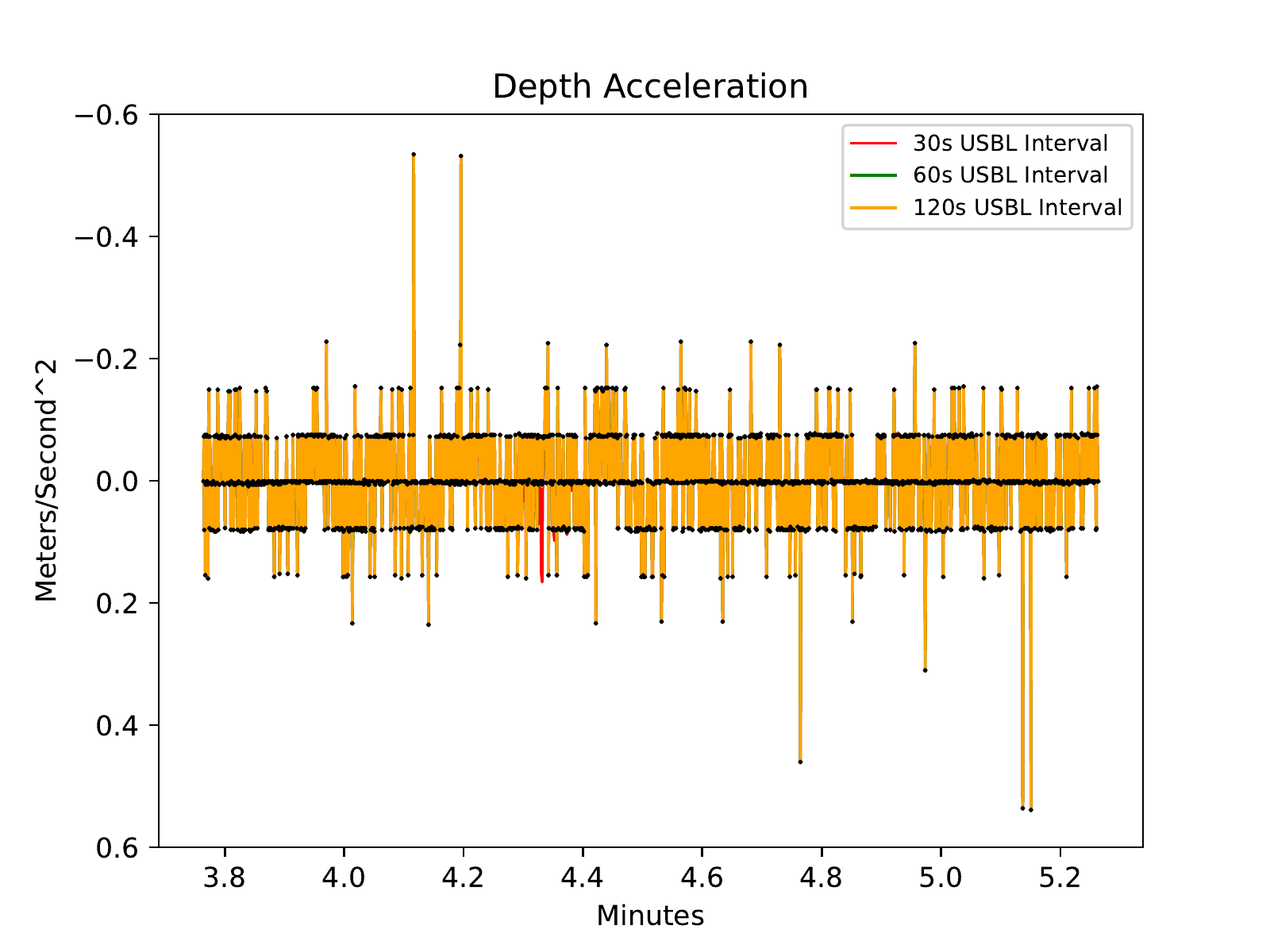}
\caption{Depth acceleration data and fit after debiasing.}
\label{fig:biasremoved}
\end{center}
\end{figure}
Figure~\ref{fig:biasremoved} shows the depth acceleration data after the bias is removed, now centered around $0$. 
Biases computed for 30, 60, and 120 second intervals were all near $0.073$.
\begin{figure}[h!]
\begin{center}
\includegraphics[scale=.5]{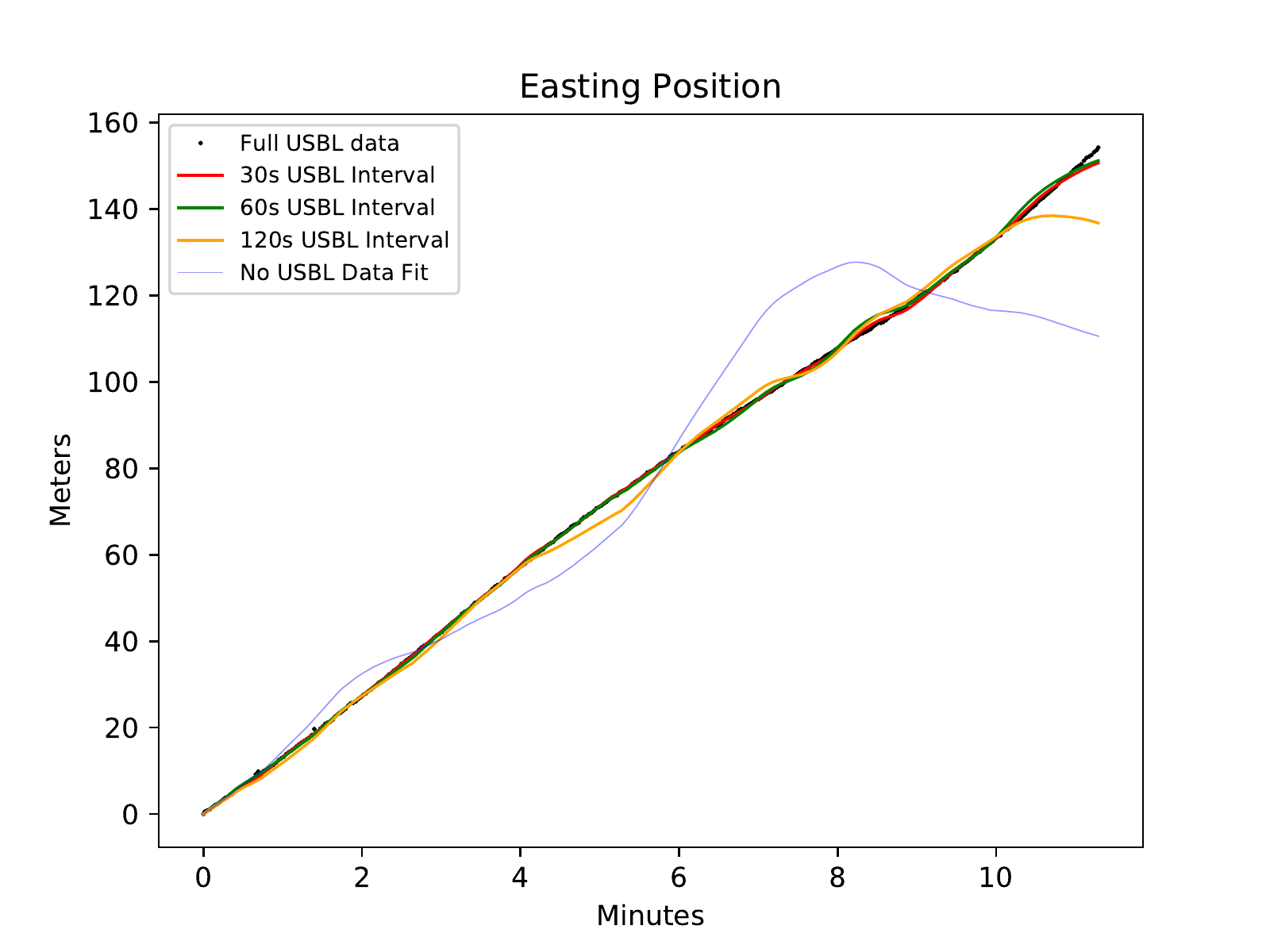}
\includegraphics[scale=.5]{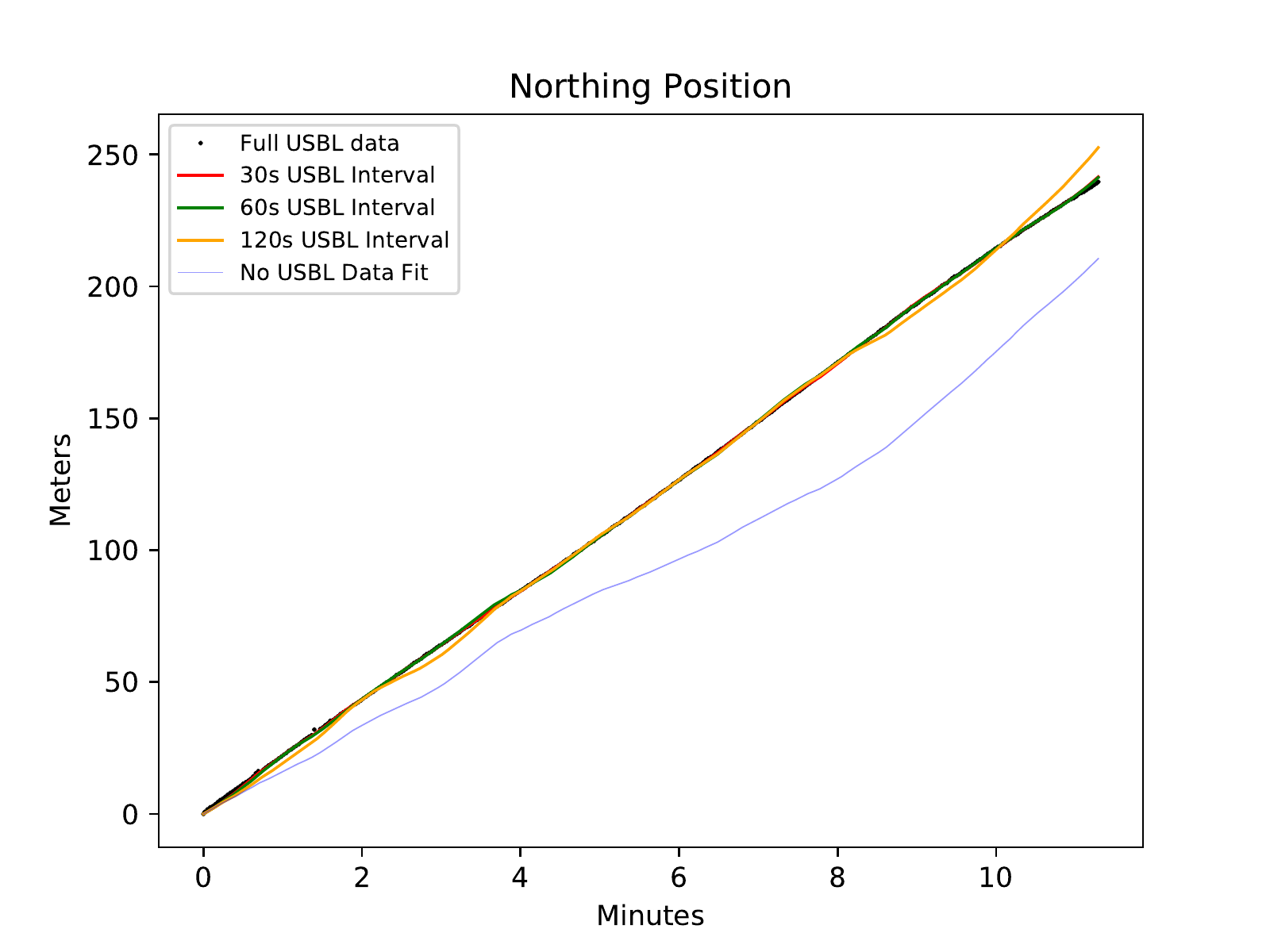}
\includegraphics[scale=.5]{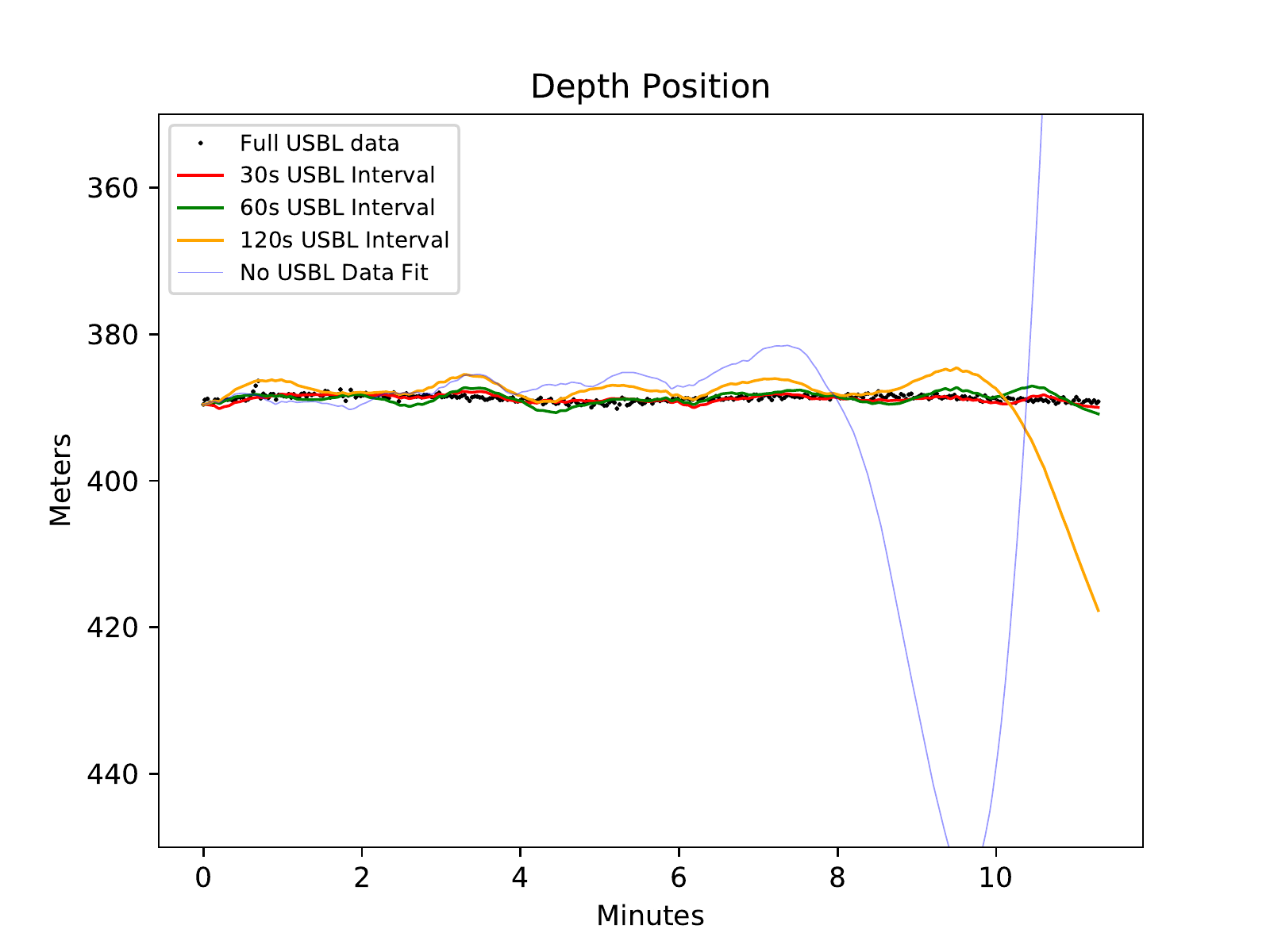}
\caption{Fitted position for three frequencies of position data. With supplemental position data the estimates perform much better than when only acceleration data is used.}
\label{fig:posshort}
\end{center}
\end{figure}
Figure~\ref{fig:posshort} has the fitted position plots for all three frequencies. The depth plot shows why using only acceleration data is can lead to large errors; small errors in acceleration data build up to have a large effect over time. However when the acceleration data is combined with a small amount of position data all three perform very well. In fact there is not a large difference in the estimates produced; this gives a promising view 
toward an online implementation.  Figure~\ref{fig:velshort} shows the fitted
velocity for all three models. Here the small differences in the fit become apparent with lower frequency position data leading to much larger changes in velocity over time. \\
\\  
{\bf Results for $50$ Minute Track}
At this scale, we consider position data at intervals of $3$ and $5$ minutes. There is also a gap in the position data near minute $27$. 
  \begin{figure}[h!]
\begin{center}
\includegraphics[scale=.5]{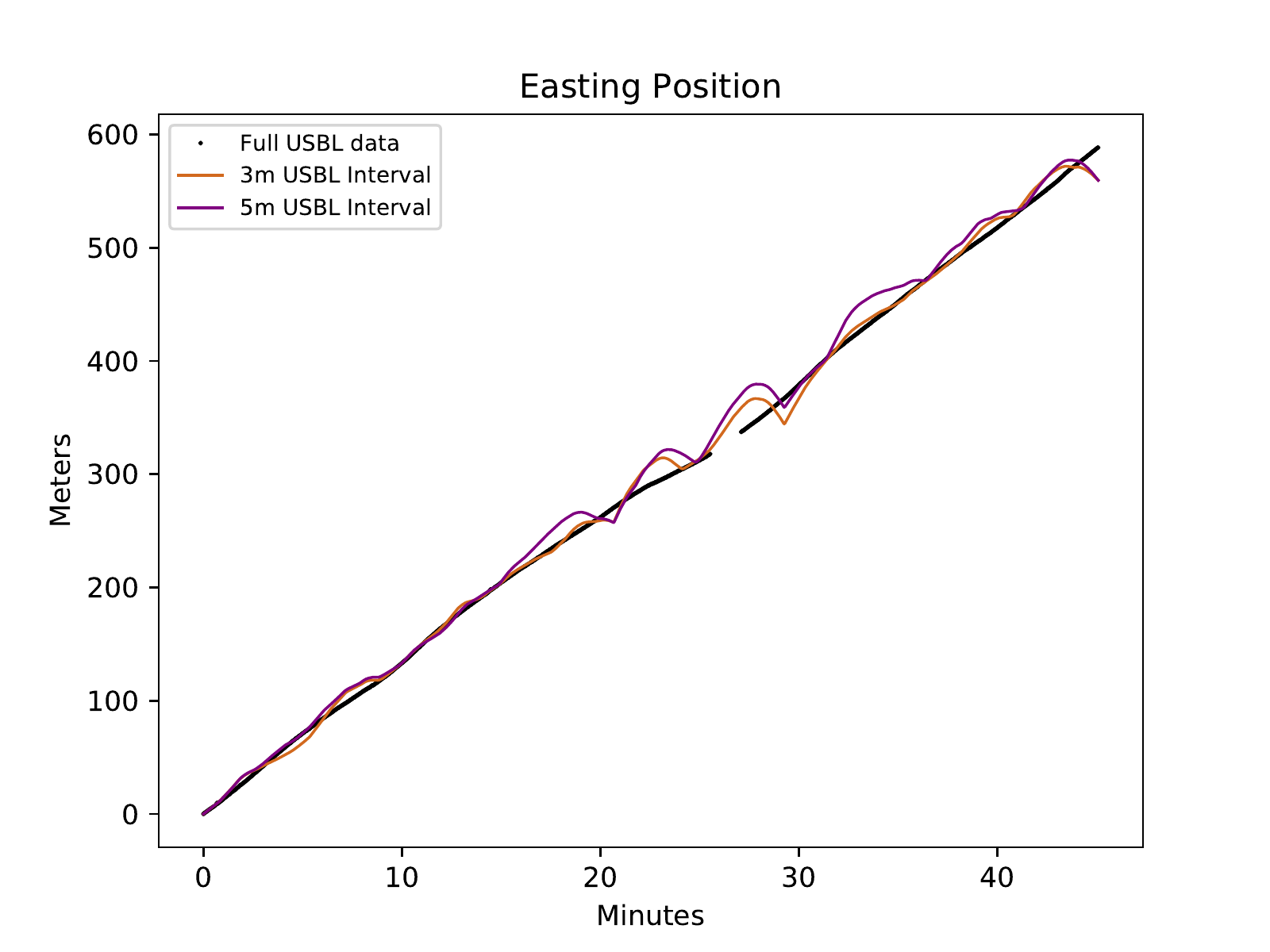}
\includegraphics[scale=.5]{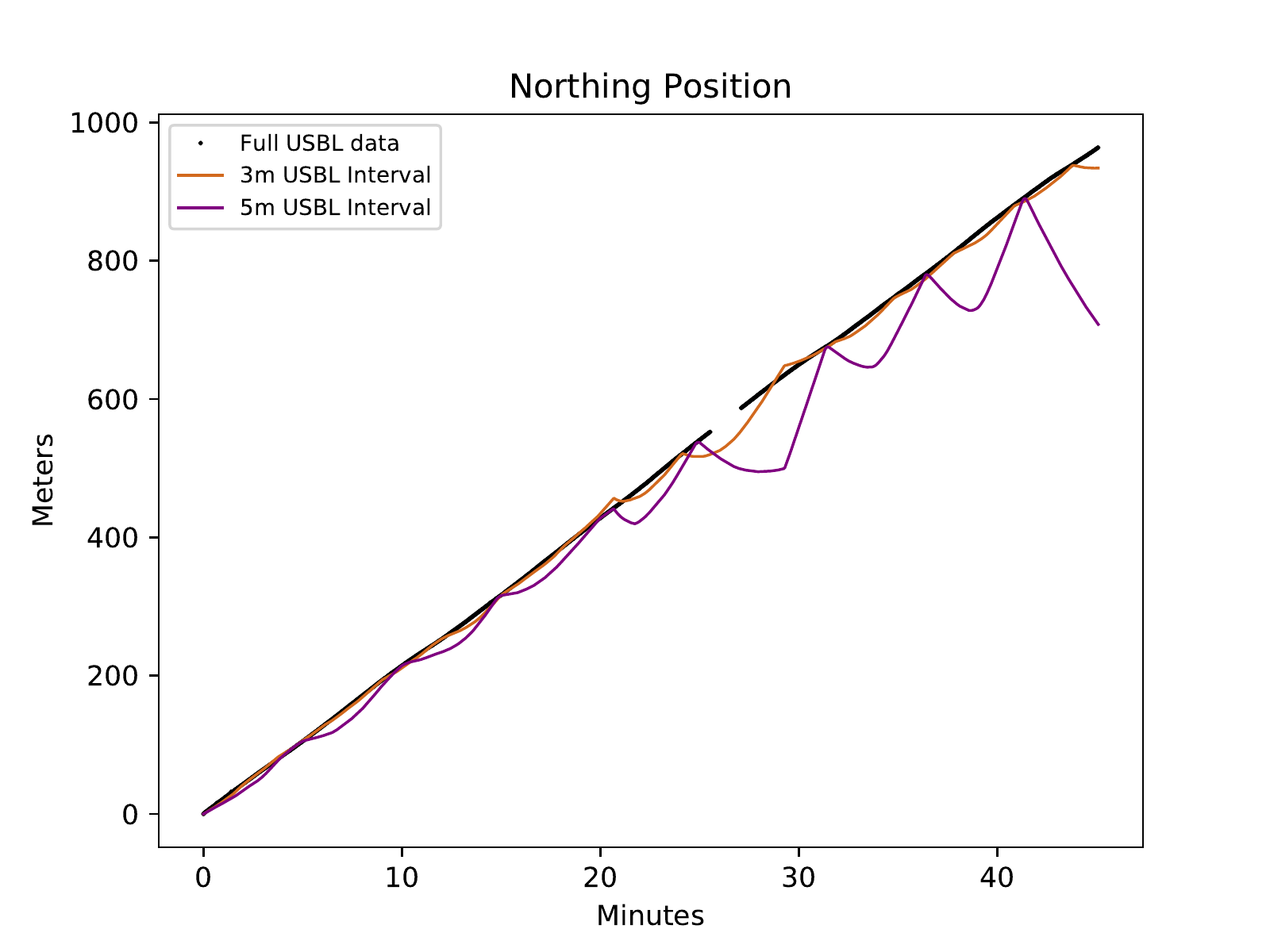}
\includegraphics[scale=.5]{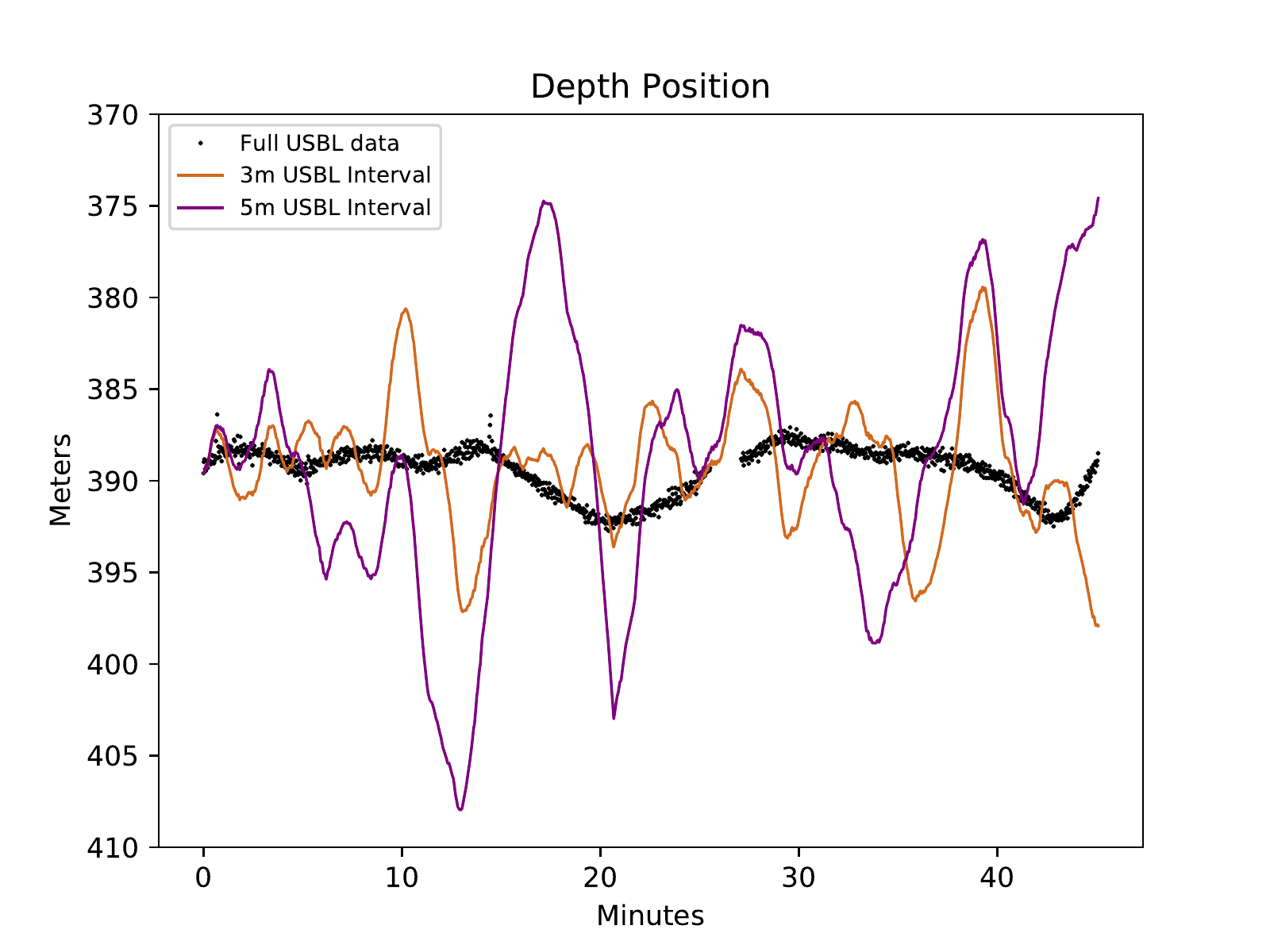}
\caption{Position fit over $50$ minutes with $3$ and $5$ minute gaps in position data. With very low frequency of position data poor acceleration in the depth data leads to much less stable position estimates.}
\label{fig:poslong}
\end{center}
\end{figure}
Figure~\ref{fig:poslong} shows the position estimates for the longer time period. When the position data is only seen every $5$ minutes the estimate becomes unstable, especially for depth, where the acceleration data quality is poorest. 
However even with three minute gaps in between position data the model performs
fairly well. Modern underwater vehicles are well-instrumented in depth, and
typically have some model for velocity (e.g. hydrodynamic velocity model in
gliders, prop counts in  propeller-drive vehicles); an extension of the methods proposed here could enable an online navigation system that requires 
ever fewer high-fidelity external position fixes (such as those provided here from the USBL data).\\

\section{Discussion}

\vspace{-.1in}

We propose a singular Kalman smoothing framework that can use 
singular covariance models for process and measurements, convex robust 
losses, and state-space constraints. The modeler can use any convex 
loss that has an implementable prox; in particular any piecewise linear-quadratic
loss and simple polyhedral constraint can be used. The framework offers a range of tools
that we illustrated using a sea survey analysis. Future work will consider real-time implementation, 
as well as extension to nonlinear models. \\
Numerical experiments illustrate that the local linear rate we have in theory requires 
a good initialization in practice. All experiments in the paper were initialized by propagating the state
estimate forward; this worked far better than an arbitrary initialization (e.g. at the $0$ vector). 
Smarter initialization can be developed for streaming/online contexts, where recent estimates play 
a key role in initializing smoothing subproblems.  

\section{Acknowledgements}
This material is based upon work supported by the Defense Advanced Research
Agency (DARPA) and Space and Naval Warfare Systems Center Pacific (SSC Pacific)
under Contract No. N66001-16-C-4001.
The work of Dr. Aravkin was supported by the Washington Research Foundation Data
Science Professorship. 


%

%
%
%

%% file: appendix.tex
\section*{Appendix}

\subsection{Proof of Theorem~\ref{thm:surjectivity}}

Conditions $2, 3, 4$ are easily seen to be equivalent. 
To see that 2 and 3 are equivalent, note that the matrix 
\[
\begin{bmatrix}
Q_i^{1/2} & 0 \\
0 & R_i^{1/2}
\end{bmatrix}
\]
is symmetric, so its nullspace is perpendicular to its range. 
Therefore surjectivity of $D_i$ is equivalent to the condition that 
the range of $\begin{bmatrix} I\\ H_i\end{bmatrix}$ covers this nullspace. \\
To see the equivalence of 2 and 4, recall that $B$ is surjective if and only if $BB^T$
is invertible, so $D_i$ is surjective exactly when the matrix 
\[
\begin{bmatrix} 
Q_i + I & H_i^T \\
H_i & R_i + H_i H_i^T
\end{bmatrix}
\]
is invertible. $Q_i +I$ is always invertible, so invertibility of the block $2\times2$ matrix is equivalent to the 
invertibility of the Schur complement $R_i + H_i\left(I-(Q_i+I)^{-1}\right)H_i^T$.\\
It remains to show that conditions 1 and 2 are equivalent. We proceed by induction on $N$.
The base case is trivial, since for $N =1$, $A = D_1$. For the inductive case, consider that 
for $N=k$ the result holds, and write the $N=k+1$ case as
\[
\begin{bmatrix}
A_k & 0 \\ 
[ 0 \quad B_k] & D_{k+1} 
\end{bmatrix} 
\begin{bmatrix} z_1 \\ z_2 \end{bmatrix}
= \begin{bmatrix} w_1 \\ w_2 \end{bmatrix}, 
\]
and assume that $A_k$ is surjective. We then know that there exists $z_1$ that satisfies $A_k z_1 = w_1$. 
The second row can now be written explicitly as 
\[
D_{k+1}z_2 = w_2 +G_{k+1}x_k, 
\]
where $x_k$ is the last component of $z_1$. Thus $A_{k+1}$ is surjective exactly when $D_{k+1}$ is, as desired.

\subsection{Proof of Lemma~\ref{lemma:property}.}
\textit{Proof:}
\textit{As $D$ is monotone we have}

\[ 
\langle \eta^* - T\eta, D\eta^* - DT\eta\rangle \geq 0
\]

\textit{as $0 \in (D+M)\eta^*$ this implies}

\[
\langle \eta^* - T\eta,-M\eta^* - DT\eta\rangle \geq 0
\]

\textit{Now} $DT\eta = DT\eta + HT\eta - HT\eta = (H-M)\eta- HT\eta$. \textit{Thus}

\[
0 \leq \langle \eta^* - T\eta,-M\eta^* + HT\eta - (H-M)\eta\rangle
\]
\[
= \langle \eta^* - T\eta,-M(\eta^* - \eta) + H(T\eta-\eta)\rangle
\]
\[
= \langle \eta^* - \eta, -M(\eta^*-\eta) + H(T\eta-\eta)\rangle
\]
\[\
 \indent+ \langle \eta-T\eta, -M(\eta^*-\eta) + H(T\eta-\eta)\rangle
\]

\textit{By definition of $M$ we have}

\[
\langle M\eta,\eta\rangle = 0
\] 

\textit{for any $\eta$. Therefore}

\[
0 \leq \langle \eta^*-\eta,H(T\eta-\eta)\rangle + \langle \eta-T\eta,-M(\eta^*-\eta)\rangle
\]
\[
+ \langle \eta-T\eta,H(T\eta-\eta)\rangle - \langle \eta-T\eta,M(T\eta-\eta)\rangle
\]
\[
= \langle \eta^* - \eta, H(T\eta-\eta) \rangle + \langle \eta-T\eta, -M(\eta^*-\eta)\rangle - ||T\eta-\eta||_{H-M}^2
\]
\[
= \langle \eta^* - \eta, H(T\eta-\eta) \rangle + \langle M(\eta-T\eta), \eta^* - \eta\rangle - ||T\eta-\eta||_{H-M}^2
\]
\[
= \langle \eta^* - \eta, (H-M)(T\eta-\eta)\rangle - ||T\eta-\eta||_{H-M}^2
\]

\subsection{Statement of \cite{AG}, Theorem 3.3.}
For a proper closed convex function $f$, the subdifferential $\partial f$ is metrically subregular at $\bar{x}$ for $\bar{y}$ with $(\bar{x}, \bar{y}) \in $ gra $\partial f$ if and only if there exists a positive constant $c$ and a neighborhood $\mathcal{U}$ of $\bar{x}$ such that
\[
f(x) \geq f(\bar{x}) + \langle \bar{y}, x-\bar{x}\rangle + cd^2(x, (\partial f)^{-1}(\bar{y})), \quad \forall x \in \mathcal{U}.
\]

\subsection{Computing with Prox Operators}

In this section, we collect the proximal operators used in the paper. 
From simple calculus, we have 
\begin{itemize}
\item 
\(
\prox_{\frac{\gamma}{2}\|\cdot\|^2}(z) = \frac{1}{1+\gamma}z.
\)
\end{itemize}
This generalizes to easily invertible least squares terms: 
\begin{itemize}
\item
\(
\prox_{\alpha \frac{1}{2}\|Ax-b\|^2}(z) = (I + \alpha A^TA)^{-1}(\alpha A^Tb+z).
\)
\end{itemize}
For $\rho(z) = \delta_{C}(z)$, we have 
\[
\prox_{\gamma \rho}(z) = \proj_{C}(z).
\]
This gives simple formulas for the following operators: 
\begin{itemize}
\item $\proj_{\gamma\mathbb{B}_2}(z) = \min(\|z\|,\gamma)\frac{z}{\|z\|}$.
\item $\proj_{\gamma\mathbb{B}_\infty}(z) = \min(\max(z,-\gamma),\gamma)$.
\item $\proj_{\mathbb{R}_+}(z) = \max(z,0)$.
\end{itemize}
We also have fast implementations for the following operators: 
\begin{itemize}
\item $\proj_{\gamma\mathbb{B}_1}(z)$, the 1-norm projection
\item $\proj_{\gamma \Delta}(z)$, the scaled simplex projection
\item $\proj_{\gamma \Delta_1}(z)$, the capped simplex projection. 
\end{itemize}

Next, the Moreau identity relates the prox operators for $f$ and $f^*$: 
\[
\prox_{\alpha f^*}(z) = z - \alpha\prox_{\alpha^{-1}}(\alpha^{-1}z) 
\]
This identity together with previous results yields the following operators: 
\begin{itemize}
\item $\prox_{\gamma\|\cdot\|_2}(z)$
\item $\prox_{\gamma\|\cdot\|_1}(z)$ 
\item $\prox_{\rho_h}(z)$, prox of hinge loss. 
\item $\prox_{\gamma \|\cdot\|_\infty}$
\end{itemize}

Often we add a simple quadratic to a penalty; the prox of the sum 
can be expressed in terms of the prox of the original penalty.

\begin{itemize}
\item
\(
\prox_{\alpha (f + \gamma/2 \|\cdot\|^2)}(x) = \prox_{\frac{\alpha}{1+2\alpha\gamma} f}\left(\frac{1}{1+2\alpha\gamma}x\right).
\)
\end{itemize}
This immediately gives the prox of the elastic net, which is the sum of the 1-norm and a simple quadratic.

Likewise, we can compute the prox of a Moreau envelope of a given penalty. 

\begin{itemize}
\item $\prox_{\gamma e_\alpha \rho}(z) = \frac{\alpha}{\gamma + \alpha} z + \frac{\gamma}{\gamma + \alpha} \prox_{(\gamma + \alpha)\rho }(z)$
\end{itemize}


This immediately gives us formulas for prox of the Huber, as well as smoothed variants of any other 
penalty in the collection.

%% file: Paper.bbl
\begin{thebibliography}{10}

\bibitem{anderson2007optimal}
B.~D. Anderson and J.~B. Moore.
\newblock {\em Optimal control: linear quadratic methods}.
\newblock Courier Corporation, 2007.

\bibitem{AndersonMoore}
B.~D.~O. Anderson and J.~B. Moore.
\newblock {\em {Optimal Filtering}}.
\newblock Prentice Hall, 1979.

\bibitem{Ansley}
C.~F. Ansley and R.~Kohn.
\newblock {A geometric derivation of the fixed interval smoothing algorithm}.
\newblock {\em Biometrika}, 69:486--487, 1982.

\bibitem{aravkin2017generalized}
A.~Aravkin, J.~V. Burke, L.~Ljung, A.~Lozano, and G.~Pillonetto.
\newblock Generalized kalman smoothing: Modeling and algorithms.
\newblock {\em Automatica}, 86:63--86, 2017.

\bibitem{JMLR:v14:aravkin13a}
A.~Y. Aravkin, J.~V. Burke, and G.~Pillonetto.
\newblock Sparse/robust estimation and kalman smoothing with nonsmooth
  log-concave densities: Modeling, computation, and theory.
\newblock {\em Journal of Machine Learning Research}, 14:2689--2728, 2013.

\bibitem{AG}
F.~Artacho and M.~Geoffroy.
\newblock Characterization of metric regularity of subdifferentials.
\newblock {\em Journal of Convex Analysis}, 15(2):365--380, 2008.

\bibitem{YAA}
Y.~Bar-Shalom, X.~R. Li, and T.~Kirubarajan.
\newblock {\em {Estimation with Applications to Tracking and Navigation}}.
\newblock John Wiley and Sons, 2001.

\bibitem{ybarshalom-2001a}
Y.~Bar-Shalom, X.~Rong~Li, and T.~Kirubarajan.
\newblock {\em Estimation with applications to tracking and navigation}.
\newblock John Wiley \& Sons, Inc., New York, 2001.

\bibitem{Bell2008}
B.~M. Bell, J.~V. Burke, and G.~Pillonetto.
\newblock An inequality constrained nonlinear {K}alman-{B}ucy smoother by
  interior point likelihood maximization.
\newblock {\em Automatica}, 45(1):25--33, Jan. 2008.

\bibitem{combettes2011proximal}
P.~L. Combettes and J.-C. Pesquet.
\newblock Proximal splitting methods in signal processing.
\newblock In {\em Fixed-point algorithms for inverse problems in science and
  engineering}, pages 185--212. Springer, 2011.

\bibitem{davis2016convergence}
D.~Davis and W.~Yin.
\newblock Convergence rate analysis of several splitting schemes.
\newblock In {\em Splitting Methods in Communication, Imaging, Science, and
  Engineering}, pages 115--163. Springer, 2016.

\bibitem{hyndman2002state}
R.~J. Hyndman, A.~B. Koehler, R.~D. Snyder, and S.~Grose.
\newblock A state space framework for automatic forecasting using exponential
  smoothing methods.
\newblock {\em International Journal of Forecasting}, 18(3):439--454, 2002.

\bibitem{Jaz}
A.~Jazwinski.
\newblock {\em {Stochastic Processes and Filtering Theory}}.
\newblock Dover Publications, Inc, 1970.

\bibitem{kalman}
R.~E. Kalman.
\newblock {A New Approach to Linear Filtering and Prediction Problems}.
\newblock {\em Transactions of the AMSE - Journal of Basic Engineering},
  82(D):35--45, 1960.

\bibitem{LFP}
P.~Latafat, N.~Freris, and P.~Patrinos.
\newblock A new randomized block-coordinate primal-dual proximal algorithm for
  distributed optimization.
\newblock {\em arXiv preprint arXiv:1706.02882}, 2017.

\bibitem{Ljung:99}
L.~Ljung.
\newblock {\em System Identification - Theory for the User}.
\newblock Prentice-Hall, Upper Saddle River, N.J., 2nd edition, 1999.

\bibitem{OGLB}
H.~Ohlsson, F.~Gustafsson, L.~Ljung, and S.~Boyd.
\newblock Smoothed state estimates under abrupt changes using sum-of-norms
  regularization.
\newblock {\em Automatica}, 48:595--605, 2012.

\bibitem{Oks}
B.~Oksendal.
\newblock {\em {Stochastic Differential Equations}}.
\newblock Springer, sixth edition, 2005.

\bibitem{Paige79}
C.~Paige.
\newblock Computer solution and perturbation analysis of generalized linear
  least squares problems.
\newblock {\em Mathematics of Computation}, 33:171--183, jan 1979.

\bibitem{RTS}
H.~E. Rauch, F.~Tung, and C.~T. Striebel.
\newblock {Maximum Likelihood estimates of linear dynamic systems}.
\newblock {\em AIAA J.}, 3(8):1145--1150, 1965.

\bibitem{RTRW}
R.~T. Rockafellar and R.~J.~B. Wets.
\newblock {\em {Variational Analysis}}, volume 317.
\newblock Springer, 1998.

\bibitem{tsay2005analysis}
R.~S. Tsay.
\newblock {\em Analysis of financial time series}, volume 543.
\newblock John Wiley \& Sons, 2005.

\end{thebibliography}
